\documentclass[graybox]{svmult}

\usepackage[english]{babel} 
\usepackage{mathptmx}       
\usepackage{helvet}         
\usepackage{courier}        
\usepackage{type1cm}        
\usepackage{amssymb,amsmath}         
\usepackage{makeidx}         
\usepackage{graphicx}        
\usepackage{multicol}        
\usepackage[bottom]{footmisc}
\usepackage{enumerate}
\usepackage{subfigure}

\DeclareMathAlphabet{\mathcal}{OMS}{cmsy}{m}{n}

\newcommand{\C}{\ensuremath{\mathbb{C}}} 
\newcommand{\N}{\ensuremath{\mathbb{N}}} 
\newcommand{\R}{\ensuremath{\mathbb{R}}} 

\newcommand{\eps}{\varepsilon}
\renewcommand{\phi}{\varphi}

\renewcommand{\E}{\ensuremath{\mathbf{E}}} 
\newcommand{\cov}{\ensuremath{\mathbf{c}\mathbf{o}\mathbf{v}}} 
\newcommand{\var}{\ensuremath{\mathbf{v}\mathbf{a}\mathbf{r}}} 
\renewcommand{\P}{\ensuremath{\mathbf{P}}} 
\newcommand{\Fc}{\ensuremath{\mathcal{F}}} 
\newcommand{\Ac}{\ensuremath{\mathcal{A}}} 
\newcommand{\ind}{\ensuremath{\mathbf{1}}} 
\newcommand{\Ns}{\ensuremath{\mathsf{N}}} 
\newcommand{\Nf}{\ensuremath{\mathfrak{N}}}

\newcommand{\Pois}{{\rm Pois}} 
\newcommand{\Exp}{{\rm Exp}}
\newcommand{\Ber}{{\rm Ber}}
\newcommand{\Binom}{{\rm Binom}}
\newcommand{\Unif}{{\rm Unif}}

\newcommand{\Tb}{\ensuremath{\mathbb{T}}}

\addto\captionsenglish{%
}

\newcommand{\Gc}{\ensuremath{\mathcal{G}}} 
\newcommand{\Hc}{\ensuremath{\mathcal{H}}}
\newcommand{\Kc}{\ensuremath{\mathcal{K}}} 

\newcommand{\Ss}{\ensuremath{\mathsf{S}}}
\newcommand{\Ds}{\ensuremath{\mathsf{D}}}
\newcommand{\Ks}{\ensuremath{\mathsf{K}}}

\newcommand{\sas}{\ensuremath{S\alpha S}}

\renewcommand{\S}{\ensuremath{\mathbb{S}}}

\newcommand{\Tt}{\tt T} 

\newcommand\argmin{{\rm argmin}\,}

\newcommand\dist{{\rm dist}\,}

\renewcommand{\dist}{{\rm dist}\,}

\newcommand{\Es}{\mathsf{E}}

\newcommand{\kommentar}[1]{}

\newcommand{\longj}{\mathop{\longrightarrow}\limits_{j \to \infty}}
\newcommand{\longs}{\mathop{\longrightarrow}\limits_{s \to t}}
\newcommand{\plongj}{\stackrel{\mathrm{\P}}{\longj}}
\newcommand{\plongs}{\stackrel{\mathrm{\P}}{\longs}}

\newcommand{\eqd}{\stackrel{d}{=}}

\newcommand\be{\begin{equation}}
\newcommand\ee{\end{equation}}
\def\qed{\hfill\hbox{${\vcenter{\vbox{
    \hrule height 0.4pt\hbox{\vrule width 0.4pt height 6pt
    \kern5pt\vrule width 0.4pt}\hrule height 0.4pt}}}$}}

\makeindex             

\newtheorem{algo}{Algorithm}

\begin{document}

\title*{Extrapolation of Stationary Random Fields}
\author{Evgeny Spodarev, Elena Shmileva and Stefan Roth}

\institute {Evgeny Spodarev \at Ulm University, Institute of Stochastics, 89069 Ulm, Germany, \\
\email{evgeny.spodarev@uni-ulm.de} \and
Elena Shmileva \at  St.Petersburg State University, Chebyshev Laboratory, St. Petersburg 199178, Russia, \\ \email{elena.shmileva@gmail.com} \and
Stefan Roth \at Ulm University, Institute of Stochastics, 89081 Ulm, Germany, \\
\email{stefan.roth@uni-ulm.de}}

\maketitle

\abstract*{We introduce basic statistical methods for the
extrapolation of stationary random fields. For square integrable
fields, we set out basics of the kriging extrapolation techniques.
For (non--Gaussian) stable fields,  which are known to be heavy
tailed, we describe  further extrapolation methods and discuss
their properties. Two of them can be considered as direct
generalizations of kriging.}

\abstract{We introduce basic statistical methods for the
extrapolation of stationary random fields. For square integrable
fields, we set out basics of the kriging extrapolation techniques.
For (non--Gaussian) stable fields,  which are known to be heavy
tailed, we describe further extrapolation methods and discuss
their properties. Two of them can be seen as direct
generalizations of kriging.}

\section{Introduction}
\label{sec:1}

In this chapter, we consider the problem of extrapolation (prediction) of random fields arising mainly in geosciences, mining, oil exploration,
hydrosciences,  insurance, etc.
It is one of the fundamental tools in geostatistics  that provides statistical inference for spatially referenced variables of interest.
Examples of such quantities are
the amount of rainfall,  concentration of minerals and vegetation, soil texture, population density, economic wealth, storm insurance claim amounts, etc.

The origins of geostatistics as a mathematical science can be
traced back to the works by B. Math\'ern (1960) \cite{Matern60},
L. Gandin (1963)  \cite{Gan63}, G. Matheron (1962-63)
\cite{Ma62,Ma63}. However, the mathematical foundations were
already laid  in the paper  by A.N.Kolmogorov  (1941) \cite{Kol41}
as well as in the book by N.Wiener (1949) \cite{Wie49}, where the
extrapolation of  stationary time series was studied, whereas
their practical application is known since 1951 due to mining
engineer D. G. Krige  \cite{Kri51}. Typical practical problems to
solve are e.g. plotting the contour concentration map of minerals
(interpolation), inference of the the mean areal precipitation and
evaluation of accuracy of the estimate from spatial measurements
(averaging or generalization), selection of
 locations of new monitoring points so that the concentration
can be evaluated with sufficient accuracy (monitoring network
design).

The remainder of this chapter is divided into three sections.
Section ~\ref{sec:2} contains preliminaries about distributional
invariance properties and dependence structure of random fields.
In Section~\ref{sec:3}, we concentrate on kriging which is a
widely used probabilistic extrapolation technique for the fields
with the finite second moment. Section~\ref{sec:4} contains recent
results on the extrapolation of heavy tailed random fields with
infinite variance, namely of stable random fields.

In Sections~\ref{sec:2} and~\ref{sec:3} we mainly follow the books
\cite{CD99,Cres91,BS12,Wack03}. Section~\ref{sec:4} is based on
the  paper \cite{KSS12}, it also contains some new results for
stable fields with the infinite first moment, see
Section~\ref{subsec:4.4}.


\section{Basics of Random Fields}
\label{sec:2}

Let $(\Omega, {\cal{F}}, \P)$ be a probability space.

\begin{definition}
A \emph{random field} \index{random field}$X$  is a random function on $(\Omega, {\cal{F}}, \P)$  indexed by points of $\mathbb{R}^d,$ $d\in \N$, i.e. $X$ is a measurable mapping $X: \Omega\times\mathbb{R}^d\to \mathbb{R}$.
\end{definition}
For an introduction into the theory of random functions see e.g. \cite[Chap. 9]{BS12}.

\subsection{Random Fields with Invariance Properties}
\label{subsec:2.1}

A random field with the finite-dimensional distributions that are invariant with respect to the action of a group $G$ of transformations of $\R^d$ is called
\emph{$G$-invariant in strict sense}\index{invariance in strict sense}. In case if this invariance is given only for the first two moments of the field which are assumed to be finite
we speak about the \emph{$G$--invariance in wide sense}\index{invariance in wide sense}. Thus, if $G$ is the group of all
translations of $\R^d$ then one calls such random fields \emph{stationary}\index{stationarity} (in respective sense). For $G$
being the group of rotations $SO_d$ one claims the random field to be \emph{isotropic}\index{isotropy}.
If $G$ is the group of all rigid motions then such field is called \emph{motion invariant}\index{motion invariance}.
The same notions of invariance can be transferred to the increments of random fields. In this case, the stationarity is often called
\emph{intrinsic}. The intrinsic stationarity in wide sense is called \emph{intrinsic stationarity of order two}\index{intrinsic stationarity of order two}. For more details on invariance properties confer \cite[Sect. 9.5]{BS12}.

\begin{exercise}
Show that the expectation (if it exists) of any process  ($d=1$) with stationary increments is a linear function, i.e.,
$\E X(t)=a\cdot t+c$ for all $t\in \R$, $a\in \R$, $c\in \R$.
\end{exercise}
A popular class of random fields are \emph{Gaussian fields}.

\begin{definition}
A random field $X=\{ X(t),\;  t\in \R^d\}$  is \emph{Gaussian} if
all its finite dimensional distributions are Gaussian.
\end{definition}

Their use for modelling purposes  in applications is explained
mainly by the simplicity of their construction and analytic
tractability combined with the normal distribution of marginals
which describes many real phenomena due to the Central Limit
Theorem.

By Kolmogorov's theorem, the probability law of a Gaussian random
field is defined uniquely by its mean value and covariance
function; see \cite[Sect. 9.2.2]{BS12} for more details. If the
mean value function $\E\, X(t),$ $t\in\R^d$ is identically zero we
call $X$ to be \emph{centered}\index{centered}. Without loss of
generality we tacitly assume all random fields of this chapter to
be centered.

\begin{exercise}
Show that for Gaussian random fields  stationarity (isotropy,
motion invariance) in strict sense and stationarity (isotropy,
motion invariance) in wide sense are equivalent. In this case we
call a Gaussian field just \emph{stationary} (\emph{isotropic},
\emph{motion invariant}).
\end{exercise}

\paragraph{\bf Examples of Gaussian Random Fields}\label{ex:WMcov}

\paragraph{\bf 1. Ornstein-Uhlenbeck Process}\index{Ornstein-Uhlenbeck process}
A centered Gaussian process $X=\{ X(t),\;  t\in \R\}$ with the
covariance function $\E \left(X(s)X(t)\right)=e^{-|s-t|/2}$,
$s,t\in\R$ is called \emph{Ornstein-Uhlenbeck Process}. Breiman
(1968) \cite[p. 350]{Brei68} has shown that $X$ is the only
stochastically continuous stationary Markov Gaussian process.
Additionally, it has short memory, i.e., $$X(t)\eqd
e^{-t/2}X(0)+V(t), \quad t>0,$$ where $V(t)$ does not depend on
the past $\{ X(s), s\leq 0 \}$, cf. \cite[Example 2.6,
p.11]{Lif12}. Defined on $\R_+$, $X$ is the strong solution of the
Langevin stochastic differential equation
$$
d X(t)=-1/2 X(t)dt+ d W(t)
$$
with initial value $X(0)\sim N(0,1)$, where $W=\{ W(t), \; t\ge
0\}$ is the standard Wiener process, see e.g. \cite[Chapt. 8,
Theorem 7]{Bulinski2}. It holds also $X \eqd \left\{
e^{-t/2}W\left( e^{\; t}\right), \; t\in\R\right\}$, cf.
\cite[Chapt. 3, p.107]{Bulinski2}.

\paragraph{\bf 2. Gaussian Linear Random Function}\index{Gaussian linear random function}
A \emph{Gaussian linear random function} $X=\{X(t), \; t\in l_2\}$
is defined by $X(t)=\langle N,t\rangle_2$, $t\in l_2$, where
$N=\{N_i\}_{i=1}^{\infty}$ is an {\it i.i.d.} sequence of
$N(0,1)$-random variables, and $l_2$ is the Hilbert space of
sequences $t=\{t_i\}_{i=1}^{\infty}$ such that
$\|t\|^2_2:=\sum_{i=1}^{\infty}t_i^2<\infty$  with scalar product
$\langle s,t\rangle_2=\sum_{i=1}^{\infty}s_i t_i$, $s,t\in l_2$.
Since $N$ is not an element of $l_2$  a.s., the expression
$\langle N,t\rangle_2$ is understood formally as the series
$\sum_{i=1}^{\infty}N_i t_i$ which converges in the mean square
sense:
$$\E\left|\sum_{i=n}^{m}N_i t_i\right|^2=\sum_{i=n}^{m}t_i ^2\rightarrow 0, \quad n,m \rightarrow \infty.$$
It holds $$X(t)\sim N(0,\|t\|^2_2), \quad X(t)-X(s)=X(t-s), \quad
\E \left(X(s) X(t)\right)=\langle s,t\rangle_2, \quad s,t\in l_2.$$
Its \emph{variogram}\index{variogram} $\gamma(h):=1/2 \cdot \E (X(t+h)- X(t))^2$ can be computed as
 $$\gamma(h)=\frac{1}{2}\E[X(h)]^2=\frac{\|h\|_2^2}{2},\quad h\in l_2,$$
see more about variograms in Sect. \ref{subsec:2.2.2}.
Here we have $\gamma(h)\rightarrow \infty$ as $\|h\|_2 \rightarrow \infty$. Transferring the notions of stationarity from
the index space $\R^d$ to $l_2$, it is clear that $X$ is intrinsic stationary of order two but not wide sense stationary.
Confer \cite{IbrRoz} for the general theory of Gaussian random functions on  Hilbert index spaces.

\paragraph{\bf 3. Fractional Brownian Field} \index{fractional Brownian
field} A \emph{fractional Brownian field} $X=\{ X(t),\;  t\in
\R^d\}$ is a centered Gaussian field with covariance (see more
about covariance in Sect. \ref{subsec:2.2.1})
$$\E (X(s)X(t))=\frac{1}{2}\left( \|s\|^{2H} + \|t\|^{2H} - \|s-t\|^{2H}\right), \quad s,t\in\R^d$$
for some $H\in(0,1]$ where $\|\cdot\| $ is the Euclidean norm in
$\R^d$. Parameter $H$ (often called \emph{Hurst index}\index{Hurst
index}) is responsible for the regularity of the paths of $X$. The
greater $H$, the smoother are the paths. For $d=1$, $X$ is called
the \emph{fractional Brownian motion}, including the {\it
two--sided} Wiener process (defined on the whole $\R$) if $H=1/2$.
In the case $d>1$, $H=1/2$ it is called the \emph{Brownian L\'evy
field}\index{Brownian L\'evy field} (see, e.g.,
\cite[Sect.~2]{Lif12}).

It is easy to check that $X$ is intrinsically stationary  of order
two and isotropic.  Its variogram
 $\gamma(h)=1/2\cdot \|h\|^{2H}$ is clearly motion invariant.
 This field is not wide sense stationary as its variance is not constant.
\begin{exercise}
 Show that $X$
\begin{enumerate}
\item has stationary increments which are positively correlated
for $H\in(1/2,1)$ and negatively correlated for $H\in(0,1/2)$.
\item is $H$--self similar, i.e., $X(\lambda t)\eqd |\lambda|^H
X(t)$ for all $\lambda\in\R$ and $t\in\R^d$. \item has a version
with a.s. H\"older continuous paths of any order $\beta\in(0,H)$.
\item has nowhere differentiable paths for any $H\in(0,1)$. \item
is a linear process for $d=H=1$, i.e., $X(t)\eqd t X_0$, $t\in\R$
for a random variable $X_0\sim N(0,1)$.
\end{enumerate}
\end{exercise}


\paragraph{\bf Examples of Non-Gaussian Random Fields}

\paragraph{\bf 1. L\'evy Process with Finite Second Moments} \index{L\'evy process}
Let $X=\{X(t),\; t\ge 0 \}$ be  a \emph{ L\'evy process} with
finite second moments. It is usually defined via the \emph{
L\'evy--Khinchin triplet} coding its jump structure, see e.g.
\cite{Sat99}. It is clear that $X$ is intrinsic stationary of
order two, but not wide sense stationary. For each of these
processes one can calculate the variance of increments and the
variogram, for example,
$$
\gamma(h)=1/2 \cdot \E(X(t+h)-X(t))^2= \lambda |h| /2,\quad h,t\ge
0
$$
for the stationary Poisson point process with intensity
$\lambda>0$.

\paragraph{\bf 2. Poisson Shot Noise Field} \index{Poisson shot noise
field} A \emph{Poisson shot noise field} $X=\{X(t),\; t\in\R^d
\}$ is defined by
$$
X(t)=\sum_{x_i\in \Phi} f(t-x_i)= \int_{\R^d} f(t-x)\Phi(dx),\ t\in \R^d,
$$
where  $\Phi$ is a  stationary Poisson point process on $\R^d$
with intensity $\lambda$, $f\in L^1(\R^d)$. It follows from
\cite[Exercise 9.10]{BS12} that $X$ is strictly stationary.

\begin{figure}[ht!]
\subfigure[Gaussian random field with Whittle-Mat\'{e}rn--type
covariance function (see Sect. \ref{subsec:2.2.1}, Example 6), $a
= 2$, $b = \nu = 1$]{ \label{fig:2}
\includegraphics[scale=.35]{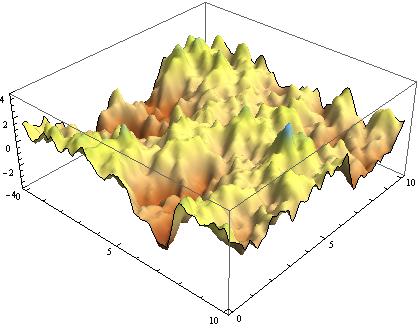}
} \hfill \subfigure[Poisson shot noise field with $\lambda = 1$
and $f(x) = \frac{1}{2\pi}\left(1-\frac{1}{4}\|x\|^2\right) \ind
(\|x\| \leq 2)$]{
\includegraphics[scale=.40]{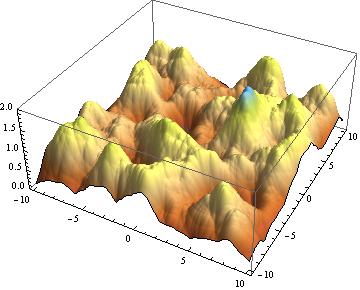}
\label{fig:1} } \caption{Simulated realizations of (strictly and
wide sense) motion invariant random fields.}
\end{figure}

It can be shown that
$$\E X(t)= \lambda \int_{\R^d} f(x) \,dx,$$
and if additionally $f\in L^2(\R^d)$ then
$$\cov \left(X(s), X(t)\right)= \lambda \int_{\R^d} f(t-s+x) f(x) \, dx,$$
i.e., the Poisson shot noise field is also wide sense stationary
(cf. \cite[Exercise 9.29]{BS12}). If $f$ is rotation invariant then $X$ is isotropic of order two. See Figure~\ref{fig:1} for a realization of $X$.

\paragraph{\bf 3. Boolean Random Function} \index{Boolean random function}
Let $\{Z_t(x),\ x \in \mathbb{R}^d\}_{t \in \mathbb{R}}$ be a
family of independent lower semi-continuous random functions with
subgraphs having almost surely compact sections 
 and $\Pi=\{(x_i, t_i)\}_{i=1}^{\infty}$ be a
Poisson point process in $\mathbb{R}^d \times \mathbb{R}$ with
intensity measure $\nu_d \otimes \theta$, where $\nu_d$ denotes
the Lebesgue measure on $\mathbb{R}^d$ and $\theta$ is a
$\sigma$-finite measure on $\mathbb{R}$. The random function
$$ Z(x) = \sup\limits_{(x_k,t_k) \in \Pi} Z_{t_k}(x-x_k), \ \ \ x \in \mathbb{R}^d$$
is called a \emph{Boolean random function}. The functions $Z_t$
are referred to as \emph{primary} functions. Boolean random
functions have been introduced by D. Jeulin for modelling rough
morphologies (\cite{JeulinAndJeulin}), see for example \cite[Sect.
7.8.1]{CD99} and references therein.

\subsection{Elements of Correlation Theory for Square Integrable Random Fields}
\label{subsec:2.2}

Let us recall the following basic concepts.

\begin{definition}
A symmetric function $f:\R^d\times \R^d\to\R$ is called
\emph{positive semi--definite}\index{positive semi--definite
function} if  for any $n\in\N$, $w_1,\ldots,w_n\in \mathbb{C}$
and any $t_1,\ldots, t_n\in \mathbb{R}^d$ it holds
$$
\sum_{i,j=1}^n  w_i \bar{w}_j f(t_i, t_j)\geq 0.
$$
\end{definition}
\begin{definition}
A symmetric function $f:\R^d\times \R^d\to\R$ is called
\emph{positive definite}\index{positive definite function} if for
 any $n\in\N$, $w_1,\ldots,w_n\in \mathbb{C}$ such that
$(w_1,\ldots,w_n)^\top \neq o\in \C^n$ and any $t_1,\ldots, t_n\in
\mathbb{R}^d$ it holds
$$
\sum_{i,j=1}^n  w_i \bar{w}_j f(t_i, t_j)> 0.
$$
\end{definition}
\begin{definition}
A symmetric function $f: \R^d\times \R^d\to \R$ is called
\emph{conditionally negative semi--definite}\index{conditionally
negative semi--definite function}
 if for any $n \in \mathbb{N}$, $w_1,\ldots,w_n\in \mathbb{C}$ such that $\sum_{i=1}^n w_i=0$  and any $t_1,\ldots, t_n\in \mathbb{R}^d$ it holds
$$\sum_{i,j=1}^n w_i \bar{w}_j f (t_i,t_j)\leq 0.$$
\end{definition}

\begin{exercise}
Prove that functions $\cos(a\cdot x),\ a\in \R$, $e^{-|x|^p},\
p\in (0,\ 2]$ are positive semi--definite, whereas  $e^{-|x|^p},\
p>2$, $|\cos x|$, $a^2+\cos^2 x,\ a\in \R$ are not.
\end{exercise}

\begin{exercise}
Find a  positive semi-definite function with discrete support.
\end{exercise}


\subsubsection{\bf Covariance function} \label{subsec:2.2.1}

\begin{definition}\label{dfn:covariance}
For a real-valued random field  $X=\{ X(t), \ t\in \R^d\}$ with
$\E X(t)^2<\infty$, $t\in\R^d$, the function $C: \R^d\times\R^d\to
\R$ given by
$$
C(s,t)=\cov(X(s), X(t))= \E \left(X(s)- \mathbb{E}
X(s)\right)\left(X(t)-\mathbb{E} X(t)\right),\quad s,t\in\R^d
$$
is called the \emph{ covariance function}\index{covariance
function}.
\end{definition}
If $X$ is wide sense stationary (motion invariant), then $C(s,t)$
depends only on $s-t$ ($\|s-t\|$, respectively), $s,t\in \R^d$.
For the properties of the covariance function see \cite[Sect.
9.4-9.6]{BS12}. We mention just a few:

\begin{enumerate}
\item \emph{Generic property}. A function $f:\R^d\times\R^d\to \R$
is a covariance function of some square integrable random field
iff it is positive semi--definite.
\begin{exercise}
Prove this fact. {\it{Hint:}} Calculate the variance of linear combinations $\sum_{i=1}^{n}x_i  X(t_i)$ for arbitrary $n\in \N$, $t_i\in \R^d$, $x_i\in \R$.
\end{exercise}

\item \emph{Spectral representation}. By Bochner-Kchinchin theorem
(see, e.g., \cite{Bochner} or \cite[Theorem 9.6]{BS12}), any
continuous at the origin positive semi--definite function  $f:
\R^d\to \R$ is a Fourier transform of some symmetric finite
measure $\mu_f$ on $\R^d$. Thus for a wide sense stationary mean
square continuous field $X$  with covariance function $C$ we have
$$\cov (X(s), X(t))=C(s-t)=\int_{\R^d} e^{i \langle x, s-t \rangle}\,\mu_C(dx).$$
Here $ \langle \cdot,\cdot\rangle$ is the Euclidean scalar product in  $\R^d$.
Measure $\mu_C$ is called a \emph{spectral measure}\index{spectral
measure} of $X$. If $\mu_C$ is absolutely continuous with respect
to the Lebesgue measure, then its density is called a
\emph{spectral density}\index{spectral density}. The above field
$X$ has itself the \emph{spectral representation}\index{spectral
representation}
\begin{equation}\label{eq:spectralRepr}
X(t)=\int_{\R^d} e^{i\langle x,t\rangle}\, \Lambda (dx),
\end{equation}
where $\Lambda(\cdot)$ is a  complex-valued orthogonal random
measure with $\E\, \Lambda(A)=0$  and $\E \left(
\Lambda(A)\overline{\Lambda(B)}\right) =\mu_C(A\cap B)$ for any
Borel sets $A,B\subset \R^d$. The integral in
\eqref{eq:spectralRepr} is understood in the mean square sense,
i.e. its integral sums converge in $L^2(\Omega,{\cal{F}}, \P)$.
For more details on the spectral representation of stationary
processes see \cite[Sect. 7, \textsection  9, \textsection
10]{Bulinski2}, \cite[Sect. 3.2, pp. 20-21]{Lif12} or
\cite[Sect.~2.3.3]{CD99}, \cite[Sect. 4.2, p. 90]{Wentzell}. The
spectral representation is used e.g. to simulate stationary
Gaussian random fields approximating the integral in
\eqref{eq:spectralRepr} by its finite integral sums with respect
to a Gaussian white noise measure $\Lambda$.
\end{enumerate}

\paragraph{\bf Parametric Families of Covariance Functions}

\paragraph{\bf 1. White Noise Model} \index{white noise}
\begin{eqnarray*}
  C(s,t)= \begin{cases}
                     \sigma^2 ,\ s=t \\
        0,\ s\neq t.
                     \end{cases}, \quad s,t\in\R^d.
 \end{eqnarray*}
It is a covariance function of a random field $X$ consisting of
independent random variables $X(t)$, $t\in\R^d$, $d\ge 1$  with
variance $\sigma^2>0$.

 \paragraph{\bf 2. Normal Scale Mixture} \emph{}\index{normal scale
mixture}
$$C(s,t)=\int_0^\infty e^{-x \|s-t \|^2} \, \mu(dx),\quad s,t\in\R^d$$
for some finite measure $\mu$ on $[0,\infty)$
 is the covariance function of a motion invariant random field for any $d\ge 1$ (see \cite{Schoenberg}).

\paragraph{\bf 3. Bessel Family} \index{Bessel family}
$$C(s,t)=b(a\|s-t \|)^{-\nu}J_{\nu}(a\|s-t \|), \quad \nu = \frac{d-2}{2},
\quad a,b>0,\quad s,t\in\R^d,$$ where $$ J_{\nu}(r) =
\sum_{j=0}^{\infty}
\frac{(-1)^j}{j!\Gamma(\nu+j+1)}\left(\frac{r}{2}\right)^{\nu+2j},\
\ r\in\R$$ is the Bessel function of the 1st kind of order $\nu$
(cf. \cite{Magnus}) and $d\ge 1$. The positive semi--definiteness
of $C$ is proven in \cite[p. 367]{Yaglom1}. The spectral density
of $C$ is given by
\begin{equation}\notag
f(h)=\frac{b(a^2-h^2)^{\nu-\frac{d}{2}}}{2^{\nu}\pi^{\frac{d}{2}}a^{2\nu}\Gamma(\nu+1-\frac{d}{2})}I(h\in[0,a]).
\end{equation}
A special case of $d=3$, i.e., $\nu=\frac{1}{2}$ yields the
so-called \emph{hole effect model}\index{hole effect model}
$$C(s,t)=b\frac{\sin(a\|s-t \|)}{a\|s-t \|},\quad s,t\in\R^d.$$
This model is valid only for $d\le 3$.

\paragraph{\bf 4. Cauchy Family} \index{Cauchy family}
\begin{equation}\notag
C(s,t)=\frac{b}{(1+(a\|s-t \|)^2)^{\nu}},\quad a,b,\nu > 0, \quad
s,t\in\R^d.
\end{equation}
Up to scaling, this function is positive semi-definite as a normal
scale mixture with $\mu(dx)=cx^{\nu-1}e^{-x}dx$ for some constant
$c>0$.

\paragraph{\bf 5. Stable Family} \index{stable family}
\begin{equation}\notag
C(s,t)=be^{-a\|s-t \|^\nu},\quad \nu \in(0,2], \quad s,t\in\R^d.
\end{equation}
This function is positive semi-definite for all $d\ge1$ since it
is made by substitution $\theta  \mapsto \|s-t \|$  out of the
characteristic function of a symmetric $\nu$-stable random
variable, cf. Definition \ref{dfn:st_var_ch_f}. A special case
($\nu=2$) of the stable family is a \emph{Gaussian
model}\index{Gaussian covariance family}: $C(s,t)=be^{-a\|s-t
\|^2}$. Its spectral density is equal to
$f(h)=\frac{b\sqrt{a}}{2}he^{-\frac{ah^2}{4}}$.

 \paragraph{\bf 6. Whittle-Mat\'{e}rn Family} \index{Whittle-Mat\'{e}rn family}
\begin{equation}\notag
C(s,t)=W_{\nu}(\|s-t\|)=b 2^{1-\nu} (a\|s-t \|)^{\nu}K_{\nu}(a \|s-t \|),
\quad s,t\in\R^d, \; s\neq t,
\end{equation}
where $\nu, a, b>0$, $d\ge 1$ and $K_{\nu}$ is the \emph{modified
Bessel function of third kind}, also called \emph{Macdonald
function}:
$$K_{\nu}(r)=\frac{\pi}{2
\sin(\pi\nu)}(e^{i\frac{\pi}{2}\nu}J_{-\nu}(re^{i\frac{\pi}{2}})-e^{-i\frac{\pi}{2}\nu}J_{\nu}(re^{-i\frac{\pi}{2}})),
\quad r\in\R, \quad \nu\not\in\N.$$ For $\nu=n\in\N$ the above
definiton of $K_{\nu}$ is understood in the sense of a limit as
$\nu\rightarrow n$, see \cite[p. 69]{Magnus}. For $s=t$, we set $C(t,t)=b.$  The spectral density
of $C$ is given by
\begin{equation}\notag
f(h)=b\frac{2\Gamma(\nu+\frac{d}{2})}{\Gamma(\frac{d}{2})\Gamma(\nu)}\frac{(ah)^{d-1}}{(1+(ah)^2)^{\nu+\frac{d}{2}}}I(h>0).
\end{equation}
If $\nu=\frac{2d+1}{2}$ then a random field with covariance
function $C$ is $d$ times differentiable in mean-square sense. If
$\nu=\frac{1}{2}$ then the \textit{exponential
model}\index{exponential model} $$C(s,t)=be^{-a\|s-t \|}, \quad
s,t\in\R^d$$ is an important special case. The same exponential
covariance belongs to the stable family for $\nu=1$.

Figure~\ref{fig:2} shows a realization of a centered Gaussian  random field $X$ with Whittle-Mat\'{e}rn
type covariance function.

\paragraph{\bf 7. Spherical Model} \index{spherical model}
is given for $1\le d\le 3$ by
$$C(s,t)=b\left(1-\frac{3}{2}\frac{\|s-t
\|}{a}+\frac{1}{2}\frac{\|s-t \|^3}{a^3}\right)I(\|s-t \| \le a),
\quad a,b>0,\quad s,t\in\R^d.$$ If $d=3$ the above formula yields
the volume of $B_{\frac{a}{2}}(0)\cap B_{\frac{a}{2}}(x_0)$, where
$x_o\in \R^3$, $\|x_0\|=\|s-t \|$. This is exactly the way how it
can be generalized to higher dimensions:
\begin{equation}\notag
C(s,t)=\nu_d\left(B_{\frac{a}{2}}(0)\cap
B_{\frac{a}{2}}(s-t)\right),\quad s,t\in \R^d,
\end{equation}
where $\nu_d$ is the Lebesgue measure.
The advantage of spherical models is that they have a compact
support.

\paragraph{\bf 8. Geometric Anisotropy} \index{geometric anisotropy}
It is easy to see that all covariance models considered above are motion invariant.
An example of a anisotropic covariance structure can be provided by rotating and stretching the argument of a motion invariant covariance model.
Let $C_0(\| h \|)$, $h\in\R^d$ be a covariance function of a motion invariant field where $C_0: \R^+\to \R^+$. For a positive definite $(d\times d)$--matrix  $Q$,
$$
C(h)=C_0(\sqrt{h^T Q h}), \quad h\in \R^d
$$
is a covariance function of some wide sense stationary anisotropic random field (see \cite[Chap. 9]{Wack03}).

\paragraph{\bf 9. Cyclone Model} \index{cyclone model}
For $d=3$, let
$$
C(x,y)=\frac{2^{3/2} det(S_x)^{1/4} det(S_y)^{1/4}}{\sqrt{det (S_x+S_y)}} W_{\nu} \left(\sqrt{(x-y)^T S_x (S_x+S_y)^{-1}S_y (x-y)}\right),
$$
where  $x,y\in \R^3$, $S_x=Id+x x^T$, $Id$ is a $(3\times 3)$--identity matrix and $W_{\nu}$ is the Whittel-Mat\'ern model.
In \cite[Theorem 5, Example~16]{Sch10}, it is shown that $C$ is a valid covariance function belonging to a more general class of covariances that mimic cyclones.
\begin{exercise}
Show $C$ is a covariance function of isotropic but not wide sense stationary random field, i.e.,
$C(x,y)=C(Rx, Ry)$ for any $R\in SO_3$, but $C(x,y)$ does not depend on $x-y$,  $x,y\in \R^3$.
\end{exercise}

For more sophisticated covariance models including spatio--temporal effects see e.g. \cite{Sch10} and references therein.
\subsubsection{\bf Variogram} \label{subsec:2.2.2}

\begin{definition}\label{dfn:variogram}
For a random field $X=\{X(t), t\in \R^d\}$  the following expression
$$
\gamma(t,s):= \frac{1}{2} \E (X(t)-X(s))^2,\quad  s,t\in \R^d
$$
is called a {\it{variogram}} \index{variogram} of $X$ whenever it  is finite for any $s,t\in \R^d$.
\end{definition}
For square integrable random fields $X$, it obviously holds
\begin{equation}\label{eq:variogram}
\gamma(s,t)=\frac{1}{2} \var X(s) + \frac{1}{2} \var X(t)-  \cov(X(t), X(s)) + \frac{1}{2} (\E X(s)-\E X(t))^2.
\end{equation}

If the field $X$ is intrinsic  stationary of order two  (motion invariant) then $\gamma(s,t)$ depends only on the difference $s-t$ ($\|s-t\|$,
respectively). With slight abuse of notation in these cases, we write $\gamma(s-t)$ and $\gamma(\|s-t\|)$ for functions  $\gamma: \R^d \to \R$ and $\gamma: \R_+\to \R$,
respectively. For a wide sense stationary random field $X$ with covariance function $C$ the relation \eqref{eq:variogram} reads
\begin{equation}\label{eq:variogram1}
\gamma(h)=C(0)-C(h),\quad h\in\R^d.
\end{equation}

\paragraph{\bf Basic Properties of Variograms}

Let $X$ be a random field with covariance function $C$ and variogram $\gamma$. The following properties hold:
\begin{enumerate}
\item $\gamma(t,t)=0$, $t\in\R^d $.

\item \emph{Symmetry}:  $\gamma(t,s)=\gamma(s,t)$, $s,t\in\R^d$.

\item \emph{Characterization of variograms}:
\begin{enumerate}[(a)]
    \item  A function $\gamma: \R^d\times \R^d\to \R_+$ is a variogram of some random field if $\gamma$
                 is conditionally negative semi--definite, see, for example, \cite[Theorem~1]{GSS01} or
                 \cite[Sect. 2.3.3, p.61]{CD99}.
                 \begin{exercise}
                        Prove that the variogram of any intrinsic stationary random field  $X$ is a conditionally
                        negative semi--definite function. \\ {\it{Hint:}} Calculate $Var (\sum_{i=1}^n \lambda_i  X(t_i))$
                        applying (\ref{eq:variogram}) with $\sum_{i=1}^n \lambda_i=0$.
                 \end{exercise}
    \item  A continuous even function $\gamma:\R^d\rightarrow\R_+$ with $\gamma(0)=0$ is a variogram of
                 a wide sense stationary random field if  $e^{-\lambda \gamma}$ is a covariance function for
                 all  $\lambda>0$, cf. \cite{Schoenberg38b}.
\end{enumerate}
\item \emph{Stability}: If $\gamma_1, \gamma_2$ are variograms
then  $\gamma=\gamma_1+\gamma_2$ is  a variogram as well. 
\begin{exercise}
Prove this fact. Show in particular that
$\gamma(h)=\gamma_1(h_i)+\gamma_2(h_j),$ where $h=(h_1,\ldots,
h_d)^\top\in \R^d$ and $\gamma_1,$ $\gamma_2$ are univariate
variograms, is a variogram.
\end{exercise}
\item \emph{Mixture}: Let $\gamma_x: \R^d\to \R_+$ be a variogram
of an intrinsic stationary (of order two) random field for any
$x\in \R$. Then the function
$$
\gamma(t)= \int_{\R} \gamma_x(t) \mu(dx), \quad t\in \R^d
$$
 is a variogram of some random field if $\mu$ is a measure on $\R$ and the above integral exists for any $t\in \R^d$, see \cite[Sect. 2.3.2, pp. 60-61]{CD99}.

\item If $X$ is wide sense stationary and $C(\infty):=
\lim_{\|h\|\rightarrow \infty}C(h)=0$, then it follows from
\eqref{eq:variogram1} that there exists the so-called
\emph{sill}\index{sill} $\gamma(\infty) := \lim_{\|h\|
\rightarrow\infty} \gamma(h)=C(0)$.

\item If $X$ is mean square continuous then $\gamma(h) \le c\|h\|^2$, $h\in\R^d$ for a constant $c>0$ and large $\|h\|$, see \cite[pp. 397-398]{Yaglom1}.

\item If $X$ is mean square differentiable then $\lim_{\|h\|\rightarrow \infty}\frac{\gamma(h)}{\|h\|^2}=0$, see \cite[pp. 136-137]{Yaglom2}.

\item Let $\gamma: \R\to \R_+$ be an even twice continuously differentiable function with $\gamma(0)=0$. Then
$\gamma$ is a variogram  iff $\gamma''$ is a covariation function, cf. \cite[Theorem~7]{GSS01}.

\end{enumerate}

\begin{exercise}
Show that for a variogram $\gamma$ the function $e^{\lambda\gamma}$ is a variogram for any $\lambda>0$.
\end{exercise}
\begin{exercise}
Let a bounded function $\gamma: \R^d\to \R_+$  be the variogram  of some intrinsic stationary of order two real valued random field $X$.
Consider $C(x,y)= \gamma(x)+\gamma(y)-\gamma(x-y)$, $x,y\in \R^d$. Show that $C$ is a covariance function of a random field $Z$ such that $Z(0)=0$ a.s.
\end{exercise}

\paragraph{\bf Parametric families of variograms}

Most parametric models for variograms of stationary random fields,
which are widely used in applications, can be constructed from the
corresponding families of covariance functions (such as those
described in Sect. \ref{subsec:2.2.1}) by applying the relation
\eqref{eq:variogram1} as well as stability  and geometric
anisotropy properties. Most models of variograms inherit their
names from the corresponding covariance models (e.g., exponential,
spherical one).
 One of few exceptions is the variogram corresponding to the white noise which is called \emph{nugget effect}\index{nugget effect}.

Stability property can be also used to create different anisotropy effects, for instance, the so-called
\emph{purely zonal anisotropy}\index{purely zonal anisotropy}. To explain this on an example, let
$\gamma(h)=a\gamma_1(h_x)+ b\gamma_2(h_y)+c \gamma_3(h_z)$, $h=(h_x,h_y,h_z)\in \R^3$, $a, b, c \geq 0$,
where $\gamma_i$ $i=1,2,3$ are variograms in dimension $d=1$. Then $\gamma$ is a variogram in dimension $d=3$ which allows for
different dependence ranges in three different axes directions.
An example of \emph{mixed anisotropy}\index{mixed anisotropy} models is
$$
\gamma(h)=\gamma_1(\|h\|)+\gamma_2\left(\sqrt{h_x^2+ h_y^2}\right)+\gamma_3(h_z),\ h=(h_x, h_y, h_z)\in \R^3.
$$
This is a mixture of 3D-isotropic variogram $\gamma_1$, 2D-isotropic (in the xy-plane) variogram $\gamma_2$ and a 1D-variogram $\gamma_3$.
Addition of a linear combination of $\gamma_2$ and $\gamma_3$ creates anisotropy in direction of z-axis.

See more about variograms in \cite[Chap. 2]{CD99}.


\subsubsection{\bf Statistical  Estimation of Covariances and Variograms}

The numerous approaches to estimate a covariance
function or a variogram are well described in the literature and therefore will not be reviewed
here. An interested reader can see e.g. \cite[Sect.~2.2]{CD99} and
\cite[Sect.~9.8]{BS12} and references therein.

\begin{example}
To illustrate the above theory, consider microscopic steel data (figure \ref{fig:steeldata}). This data is
obviously isotropic. Figure \ref{fig:steelvario} shows estimates for the corresponding variogram.
For this purpose Math\'{e}ron's estimator (see \cite[p. 325]{BS12}) was calculated for different directions and $0 \leq h \leq 0.5$.
The directions can be distinguished by the color of their plots. Since these data are isotropic the estimates
differ not too much.
\begin{figure}[ht!]
\subfigure[Microscopic image of a steel surface.]{ \label{fig:steeldata}
\includegraphics[scale=.40]{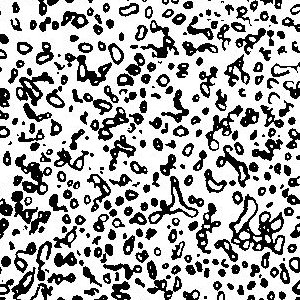}
}
\hfill
\subfigure[Estimates for the x-direction (red), y-direction (green), all directions (black) for
values $0 \leq h \leq 0.5$.]{ \label{fig:steelvario}
\includegraphics[scale=.27]{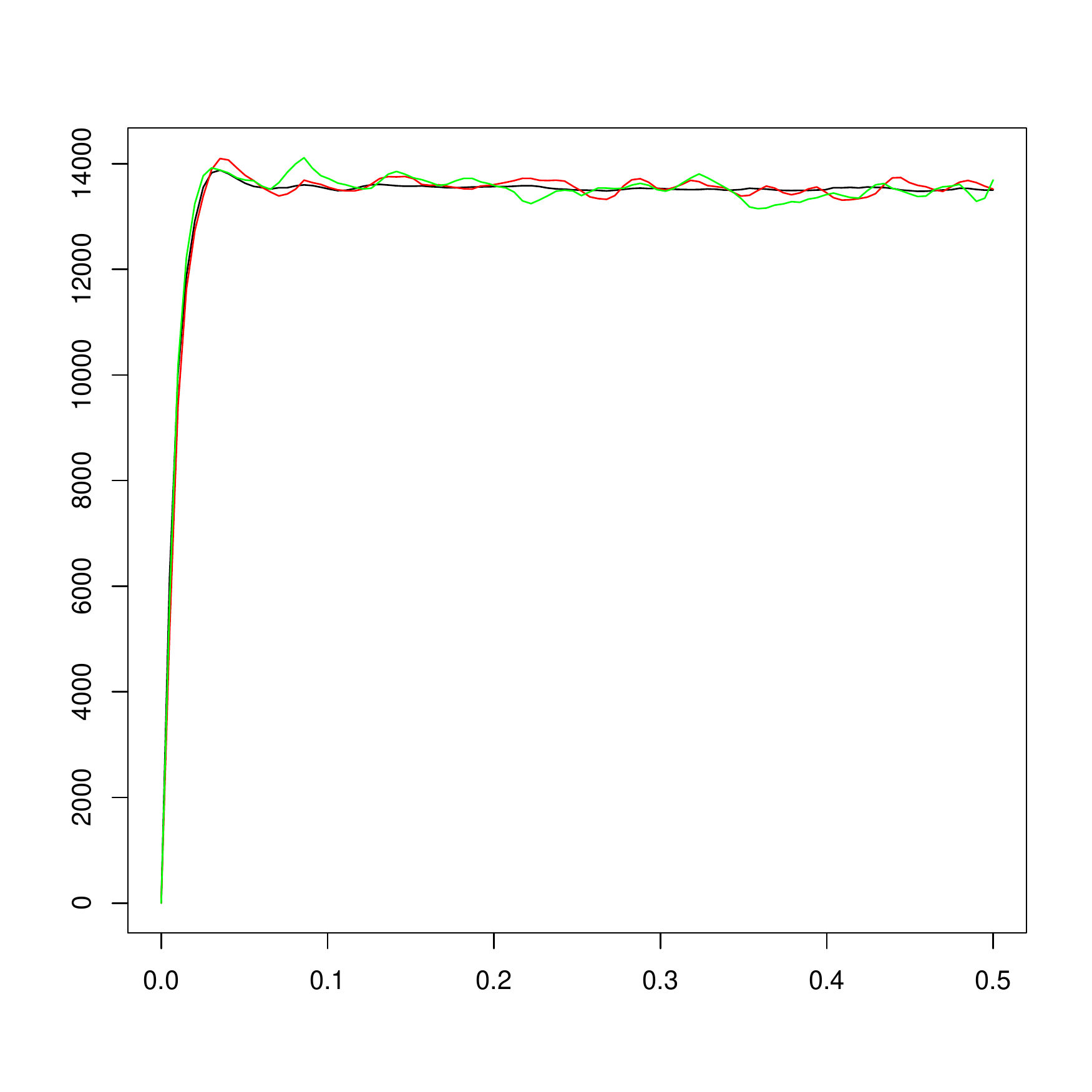}
}
\caption{Microscopic steel image (left) and its empirical variogram estimated in different directions (right)}
\label{fig:3}
\end{figure}
\end{example}

\begin{example}
Let us construct an example of zonally anisotropic variogram, in which the value for the sill depends on the direction of the input vector $h$. Consider
\begin{equation*}
    \gamma(h) = \gamma_1(h) + \gamma_2(h)
\end{equation*}
where $\gamma_1$ is an isotropic variogram
\begin{equation*}
    \gamma_1(h) = 1-e^{-|h|}, \ \ \ h \in \mathbb{R}^2
\end{equation*}
and $\gamma_2$ is a geometrical anisotropic variogram
model
\begin{equation*}
    \gamma_2(h) = 1-e^{-\frac{\sqrt{h^T Q h} }{5}}, \ \ \ h \in \mathbb{R}^2
\end{equation*}
with $Q = \sqrt{\Lambda} \cdot R$ with
$R$ being a rotation matrix with rotation angle $\alpha = 2$ and $\Lambda = diag(5,1)$
being a diagonal matrix.
Figure~\ref{fig:zonal} shows $\gamma$ on $[-1,1]^2$.
Figure \ref{fig:zonal2} illustrates the elliptic form of the contour lines of a zonally
anisotropic variogram.
\begin{figure}[ht!]
\subfigure[Zonally anisotropic model with rotation angle $\alpha = 2$ and scaling factors
$\lambda_1 = 5$, $\lambda_2 = 1$.]{ \label{fig:zonal}
\includegraphics[scale=.35]{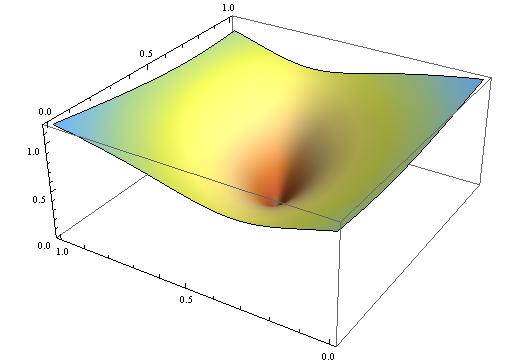}
}
\hfill
\subfigure[Contour lines of \ref{fig:zonal}]{\label{fig:zonal2}
\includegraphics[scale=.30]{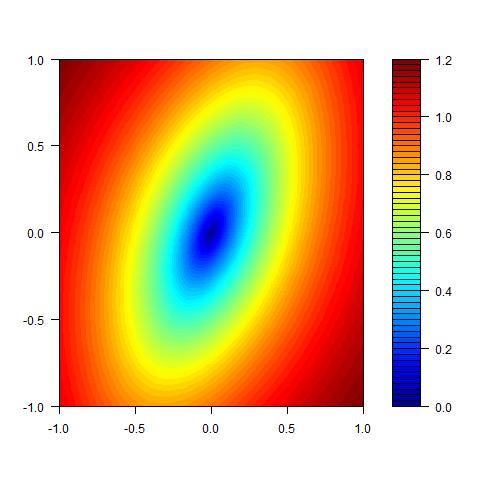}
}
\caption{Zonally anisotropic theoretical variogram}
\label{fig:4}
\end{figure}
\end{example}


\subsection{Stable Random Fields}
\label{subsec:2.3}

In this Section, we review the basic notions of the theory of
stable distributions, random measures and fields. A very good
reference which covers most of this topic is \cite{ST94}, see also
\cite{Nol}, \cite[Chapter~3]{Sat99}, \cite{Zol86}, etc.


\subsubsection{Stable Distributions}
Let $n\in \mathbb{N}$. We begin with the definition of stability for random vectors.

\paragraph{\bf Stable Random Vectors}

\begin{definition}
\label{dfn:st_vect_stab}
A random vector $\boldsymbol{X} = (X_1,\ldots,X_n)^\mathsf{T}$ in $\mathbb{R}^n$ is called  {\it{stable}}\index{stable random vector}  if
 for all $m\geq 2$  there exist $c=c(m)>0$ and $k=k(m)\in\mathbb{R}^n$ such that
$$\boldsymbol{X}^{(1)}+\boldsymbol{X}^{(2)}+...+\boldsymbol{X}^{(m)}\buildrel{d}\over{=} c \boldsymbol{X}+ k,$$
where $\{\boldsymbol{X}^{(i)}\}_{i=1}^m$ are independent copies of $\boldsymbol{X}$.
\end{definition}
It can be shown that $c=m^{1/\alpha}$ for some $0 < \alpha \leq 2$ which is called the \emph{stability index}\index{stability index}, see \cite[Theorem~2.1.2]{ST94}.
There is an equivalent  definition of stable vectors which is often used in mathematical practice to check stability.

\begin{definition}
\label{dfn:st_vect_ch_f}
 Let $\alpha\in (0,2)$. We say that a random vector $\boldsymbol{X} = (X_1,\ldots,X_n)^\mathsf{T}$ in $\mathbb{R}^n$ is  {\it{$\alpha$-stable}}
if its characteristic function is given by
\begin{eqnarray}
\label{eq:st_vector}
 \varphi_{\boldsymbol{X}}(\boldsymbol{\theta}) = \begin{cases}
                    e^{-\int_{\S^{n-1}} |\langle\boldsymbol{\theta},\boldsymbol{s}\rangle|^\alpha\left(1-i\text{sign}(\langle\boldsymbol{\theta},\boldsymbol{s}\rangle)\tan \frac{\pi\alpha}{2}\right)\Gamma(d\boldsymbol{s}) + i\langle\boldsymbol{\theta},\boldsymbol{\mu}\rangle}, &\alpha \neq 1, \\
            e^{-\int_{\S^{n-1}} |\langle\boldsymbol{\theta},\boldsymbol{s}\rangle|\left(1+i\frac{2}{\pi}\text{sign}(\langle\boldsymbol{\theta},\boldsymbol{s}\rangle)\ln|\langle\boldsymbol{\theta},\boldsymbol{s}\rangle|\right)\Gamma(d\boldsymbol{s}) + i\langle\boldsymbol{\theta},\boldsymbol{\mu}\rangle}, &\alpha = 1,
                   \end{cases}
\end{eqnarray}
where  $\Gamma$ is a  finite measure on the unit sphere $\S^{n-1}$ of $\mathbb{R}^n$ and  $\boldsymbol{\mu}$ is an arbitrary vector  in $\mathbb{R}^n$.
\end{definition}

The pair $(\boldsymbol{\mu}, \Gamma)$ gives a unique
parametrization of the distribution of $\alpha$-stable random
vectors for $\alpha\in (0,2)$, and we write $X\sim
S_{\alpha}(\boldsymbol{\mu}, \Gamma).$ This means that there is no
other pair $(\mu', \Gamma')$ yielding the same  characteristic
function $\varphi_{X}$ in (\ref{eq:st_vector}). The measure
$\Gamma$  is called \textit{spectral measure}\index{spectral
measure! of a stable vector} of $X$ and contains all the
information about the dependence between the vector components
$X_i$
 (see also Exercise~\ref{exe:srv}).
The vector $\boldsymbol{\mu}$ reflects the \textit{shift} with respect to the origin\index{shift of a stable vector}.

\begin{definition}\label{dfn:singular}
A random vector $X=(X_1,\ldots,X_n)^\top$ is called \emph{singular}\index{singular random vector} if $\sum_{i=1}^n c_i X_i = 0$
a.s. for some $(c_1,\ldots,c_n)^\mathsf{T} \in \R^n \setminus \{0\}$. Otherwise,  it is called \emph{full-dimensional}\index{full-dimensional random vector}.
\end{definition}

If $\alpha=2$, then  Definition~\ref{dfn:st_vect_stab} yields a \emph{Gaussian} random vector
which is equivalently defined  via  its characteristic function
\begin{equation}\label{eq:gaus_vect}
\varphi_{\boldsymbol{X}}(\mathbf{\theta}) = \exp\left\{ i\langle\mathbf{\theta},\boldsymbol\mu\rangle - \tfrac{1}{2} \mathbf{\theta}^\mathsf{T}\boldsymbol\Sigma \mathbf{\theta} \right\}.
\end{equation}
Here  $\mu \in \mathbb{R}^n$ is the \emph{mean} of $X$ and
$\Sigma$ is a symmetric, positive semi--definite $(n\times
n)$--\emph{covariance matrix} of $X$. The matrix $\Sigma$ has the
elements $\sigma_{{i j}}= \E (X_i-\mu_i)(X_j-\mu_j)$, where $X_i$
and $\mu_i$ are the components of vectors $X$ and $\mu$,
respectively. It is easy to see that if  $ \det \Sigma=0$ then the
Gaussian random vector $\boldsymbol{X}$ is singular\index{singular
random vector!Gaussian}.

\begin{exercise}
Prove that  Definition \ref{dfn:st_vect_stab} is equivalent to Definition~\ref{dfn:st_vect_ch_f} for $\alpha\in (0,2)$, and it is equivalent to the definition of a Gaussian random vector
via relation \eqref{eq:gaus_vect} for $\alpha=2$.
\end{exercise}

\begin{exercise}
Show that for $X\sim S_{\alpha}(\mu, \Gamma)$ the relation between
the drift $k$ in Definition~\ref{dfn:st_vect_stab} and the shift
$\mu$ in Definition~\ref{dfn:st_vect_ch_f} is
$k(m)=\mu(m-m^{1/\alpha})$. {\textit{Hint:}} First show that
$\sum_{i=1}^m X^{(i)}\sim S(m\mu, m\Gamma)$ and
$m^{1/\alpha}X+k(m)\sim S(m^{1/\alpha}\mu+k(m), m\Gamma)$.
\end{exercise}

\begin{remark}
For $\alpha=2$, the characteristic function \eqref{eq:st_vector} has the form
\begin{equation}\label{eq:gaus_vect1}
\varphi(\boldsymbol{\theta}) =  \exp\left\{-\int_{\mathbb{S}^{n-1}} \langle\boldsymbol{\theta},\boldsymbol{s}\rangle^2\Gamma(d\boldsymbol{s}) + i\langle\boldsymbol{\theta},\boldsymbol{\mu}\rangle\right\}.
\end{equation}

It is easy to find two different  finite measures $\Gamma_1$ and $\Gamma_2$  on $\S^{n-1}$  yielding the same function $\varphi$ in this case.
\end{remark}
\begin{exercise}
Check that the following two finite measures on the unite sphere in $\mathbb{R}^2$
$$\Gamma_1(ds)=  \delta_{(\sqrt{2}/2, \sqrt{2}/2)}(ds)+\delta_{(-\sqrt{2}/2, -\sqrt{2}/2)}(ds),$$
$$\Gamma_2(ds)= 2 \delta_{(\sqrt{2}/2, \sqrt{2}/2)}(ds)$$
and a shift $\mu\in \mathbb{R}^2$ yield the same expression in
(\ref{eq:st_vector}) if $\alpha=2$, $n=2$. Here
$\delta_{x}(\cdot)$ is the Dirac measure concentrated at the point
$x\in\R^2$. Verify that this expression corresponds to the
characteristic function of the Gaussian vector with shift $\mu$
and covariance matrix
\begin{equation*}
\Sigma= \begin{pmatrix}
 2 & 2 \\
 2 &  2\end{pmatrix}.
\end{equation*}
\end{exercise}

A random vector $\boldsymbol{X}$ in $\R^n$ is called \emph{symmetric}\index{symmetric distribution of a random vector} if
 $\mathbb{P}(\boldsymbol{X} \in A) = \mathbb{P}(-\boldsymbol{X} \in A)$ for any Borel set $A\in \R^n$. For symmetric $\alpha$-stable distributions,
we use the standard abbreviation  $S\alpha S$.
\begin{lemma}[\cite{ST94}, Theorem 2.4.3]
An $\alpha$--stable random vector $X$  is symmetric iff its shift
$\mu=0$ and spectral measure $\Gamma$  is symmetric.
\end{lemma}

\begin{exercise} \label{exe:srv}  Let $X=(X_1, X_2)^T$ be an $\alpha$-stable random vector, $\alpha\in (0,2)$, with the spectral measure $\Gamma$. Let $supp{(\Gamma)}$ be the support
of $\Gamma$.  Show that
\begin{itemize}
\item
 $X_1$ is independent of $X_2$ iff $supp{(\Gamma)}$ lies within the intersection of the sphere with the coordinate axes.
\item
 $X_1= c\cdot X_2$ a.s.  for some $c\in\R$ (i.e. the vector $X$ is singular) iff $supp{(\Gamma)}$ is a subset of the unite sphere intersected by a hyperplane.
\end{itemize}
\end{exercise}

\paragraph{\bf  Stable Random Variables}

If $n=1$ we deal with stable random variables whose distribution laws are defined by
four parameters $\alpha$, $\sigma$, $\beta$, and  $\mu$.
\begin{definition}\label{dfn:st_var_ch_f}
The random variable $X$ is called \emph{$\alpha$-stable}\index{stable random variable} if its
characteristic function has the form
 \begin{eqnarray*}
  \varphi_X(\theta) = \begin{cases}
                      \exp\left\{-\sigma^\alpha |\theta|^\alpha \left(1-i\beta(\text{sign}(\theta))\tan \frac{\pi \alpha}{2}\right) + i \mu\theta\right\}, & \alpha\in (0,2], \alpha \neq 1, \\
              \exp\left\{-\sigma |\theta| \left(1+i\beta \frac{2}{\pi}(\text{sign}(\theta))\ln|\theta|\right) + i \mu\theta\right\}, & \alpha = 1.
                     \end{cases}
 \end{eqnarray*}
We write $X\sim S_{\alpha}(\sigma, \beta, \mu)$.
\end{definition}
Compared with representation (\ref{eq:st_vector}), two new parameters $\sigma\geq 0$ and $\beta\in [-1,1]$
 introduced in lieu of the spectral measure $\Gamma$ are interpreted as parameters of scale and skewness, respectively.
\begin{exercise}
Show that the spectral measure of $X\sim S_{\alpha}(\sigma, \beta, \mu)$ is given by
$$
\Gamma(ds)=\frac{\sigma^{\alpha}}{2}(1+\beta)\delta_{1}(ds)+\frac{\sigma^{\alpha}}{2}(1-\beta)\delta_{-1}(ds).
$$
Hence, it holds
$$
\sigma^{\alpha}=\Gamma(\{1\})+\Gamma(\{-1\})
, \qquad
\beta=\frac{\Gamma(\{1\})-\Gamma(\{-1\})}{\Gamma(\{1\})+\Gamma(\{-1\})}.
$$
\end{exercise}
\begin{remark}
Stable distributions are absolutely continuous. Nevertheless, their densities are not known in the closed form except for the cases  $\alpha=1/2$, $\alpha=1$ and $\alpha=2$.
\end{remark}

\begin{example} ~ 
\begin{enumerate}
\item $X\sim S_2(\sigma,0,\mu)$ is a Gaussian random variable with
mean $\mu$ and variance $2\sigma^2$. \item Random variable $X\sim
S_{\alpha}(\sigma, \pm 1,\mu)$ is called  \emph{totally
skewed}\index{stable random variable!totally skewed}. Notice that
if $\alpha\in [1,2)$ then $X$ attains values in the whole $\R$. On
the contrary, if $\alpha\in (0,1)$ and $\mu=0$, then  $X\ge 0$,
$(X\le 0)$ a.s. when $\beta=1$  $(\beta=-1)$, respectively.
\end{enumerate}
\end{example}
\begin{exercise}
Show that  the characteristic function of $S\alpha S$ random variable $X$ is equal to
$\varphi_{X}({\theta}) = \exp\{-\sigma^{\alpha} |\theta|^{\alpha}\}$, i.e.,  $X\sim S_{\alpha}(\sigma, 0,0)$ for some $\sigma>0$.

\end{exercise}
\paragraph{\bf  Tails and Moments}

The non--Gaussian stable distributions are fat tailed. This means
that they belong to a subclass of heavy tailed distributions with
especially slow large deviation behavior, see more details on
heavy tailed distributions e.g. in \cite{FKZ11}, \cite{Mar07},
etc. Namely, for $X\sim S_{\alpha}(\sigma, \beta,\mu)  $  with
$\alpha\in (0,2)$ there exists $c>0$ such that
\begin{equation}\label{eq: st_tail}
\P(|X|>x)\sim c x^{-\alpha},\ x\to\infty.
\end{equation}
Here and in what follows we say that $a_x \sim b_x$ if $\lim_{x\to \infty}\frac{a_x}{b_x}=1$.
As a corollary of \eqref{eq: st_tail}, the absolute moments of $X$ behave like
 $$\E |X|^p=\int_{0}^{\infty}\P\{|X|>x^{1/p}\} dx\approx c_1 \int_0^{\infty} x^{-\alpha/p} dx.$$
They are finite if $p\in (0,\alpha)$ and
 infinite for any $p\in [\alpha,\infty)$.

\begin{exercise} \label{exe: ind} Show that
\begin{itemize}
\item normal distribution $X\sim N(\mu,\sigma^2)$ is not heavy
tailed (this is equivalent to the statement that the tails are
exponentially bounded), i.e.,
$$\P(X<-x)=\P(X>x)\sim \frac{1}{\sqrt{2\pi}\sigma x}e^{-x^2/(2\sigma^2)},\quad x\to\infty.$$
\item
for $X\sim S_{\alpha}(\sigma, \beta, 0)$, $\alpha \in (0,2)$, $\alpha \neq 1$ it holds
\begin{equation}
\left(\E\vert X \vert^p\right)^{1/p} = c_{\alpha,\beta}(p) \sigma
\label{eq:p_mean}
\end{equation}
for every $p\in (0,\alpha)$. Here
 $c_{\alpha,\beta}(p)= (\E |\xi|^p)^{1/p}$ with  $\xi\sim S_{\alpha}(1,\beta,0)$. If $\alpha=1$ then equation (\ref{eq:p_mean}) holds  only for $\beta=0$.
\item for any $\alpha$-stable random variables $X$ and $Y$  the
sum $aX+bY$, $a, b \in \mathbb{R}$ is again  $\alpha$-stable.
Moreover, components $X_i$ of the stable vector $X=(X_1, X_2) \sim
S_{\alpha}(\mu, \Gamma)$ are stable, and it holds $\sigma_{a X_1+b
X_2}=\int_{\S^1}|a s_1 +b s_2|^{\alpha} \Gamma(d s_1, d s_2)$ for
any $a,b \in \mathbb{R}$.
\end{itemize}
\end{exercise}
Simulation of stable random variables is extensively described in  \cite{Nol}.

\subsubsection{Integration with Respect to Stable Random Measures}

Let $(E,\mathcal{E},m)$ be an arbitrary measurable space with $\sigma$-finite measure $m$ and  $\mathcal{E}_0 := \{A \in \mathcal{E}: m(A) < \infty\}$. Let $\beta: E \to [-1,1]$ be a measurable function.
\begin{definition}
\label{dfn:rm}
A random function $M=\{M(A), \; A\in \mathcal{E}_0\} $ is called an
\emph{independently scattered random measure}\index{independently scattered random measure} (\emph{random noise}) if
\begin{enumerate}
\item for any $n \in \mathbb{N}$ and pairwise disjoint sets $A_1,A_2, \ldots, A_n\in \mathcal{E}_0$ random variables $M(A_1),\ldots,M(A_n)$ are independent,
\item $M(\bigcup_{j=1}^\infty A_j) = \sum_{j=1}^\infty M(A_j)$ a.s.  for a sequence of disjoint sets $A_1,A_2,\ldots \in \mathcal{E}_0$ with $\bigcup_{j=1}^\infty A_j \in \mathcal{E}_0$.
\end{enumerate}
\end{definition}
\begin{definition}
\label{dfn:alphasrm}
An independently scattered random measure $M$ on $(E, \mathcal{E}_0)$ is called \emph{$\alpha$-stable}\index{stable random measure}  if
  for each $A \in \mathcal{E}_0$
$$M(A) \sim S_\alpha\left((m(A))^{1/\alpha}, \frac{\int_A\beta(x)m(dx)}{m(A)}, 0\right).$$
Measure $m$ is called \emph{control measure}\index{stable random
measure!control measure}, and $\beta$ is the \emph{skewness
function}\index{stable random measure!skewness function} of $M$.
\end{definition}

Our goal is to define an integral $\int_E f(x) M(dx)$ of a deterministic function $f:E\to \R$ with respect to an  $\alpha$-stable random measure $M$.
For a simple function
$f(x)=\sum_{i=1}^n c_i 1_{A_i}(x)$, where $\{A_i\}_{i=1}^n\subset \mathcal{E}_0$  are pairwise disjoint, we set
$$
 \int_E f(x) M(dx)=\sum_{i=1}^n c_i  M(A_i).
$$
It can be shown that, so defined, the integral $\int_E f(x) M(dx)$ does not depend on the  representation of $f$ as a simple function, see \cite[Sect.3.4]{ST94}.
For an arbitrary $f:E \to \mathbb{R}$ such that $\int_E \vert f(x) \vert^\alpha m(dx) < \infty$  consider a pointwise  approximation of $f$
by simple functions $f^{(n)}$. Then we set
$$
 \int_E f(x) M(dx)=\mbox{plim}_{n\to \infty} \int_E f^{(n)}(x) M(dx).
$$
Here  $\mbox{plim}$ denotes the limit in probability. This
definition is independent of the choice of the approximating
sequence $\{f^{(n)}\}$, cf. \cite[Sect.~3.4]{ST94} for more
details.
\begin{lemma}
Let $X=\int_E f(x) M(dx)$, where $M$ is an $\alpha$-stable random
measure  with control measure $m$ and skewness function $\beta$.
Then $X$ is an $\alpha$-stable random variable with zero shift,
scale parameter
\begin{equation}\label{eq:scale}
\sigma_{X}^\alpha = \int_E \vert f(x) \vert^\alpha m(dx),
\end{equation}
and skewness parameter
$$
\beta_{X} =\frac{\int_E f(x) ^{<\alpha>}\beta(x) \, m(dx)} {\int_E \vert f(x) \vert^\alpha \,m(dx)},
$$
where $a^{<p>}= sign(a)\cdot |a|^p.$
\end{lemma}
For the proof see \cite[Sect.3.4]{ST94}. Notice that if
$\beta(x)=0$ for all $x \in E$ then   the integral $X$ is a
$S\alpha S$ random variable.


In case of stable vectors with an integral representation, we have
the following criterion of their full--dimensionality /
singularity.
\begin{lemma}\label{lemma:great_subsphere}
Consider a $n$--dimensional $\alpha$--stable random vector
$\vec{X}=(X_1,\ldots,X_n)^\mathsf{T}$ with $0<\alpha\leq 2$ and
integral representation
$$\vec{X} = \left(\int_E f_1(x) M(dx),\ldots,\int_E f_n(x) M(dx)\right)^\mathsf{T}.$$
Then $\vec{X}$ is singular if and only if $\sum_{i=1}^n c_i f_i(x)
= 0$ $m$--almost everywhere for some vector
$(c_1,\ldots,c_n)^\mathsf{T} \in \R^n \setminus \{0\}$.
\end{lemma}
The proof of Lemma~\ref{lemma:great_subsphere} follows from
Definition~\ref{dfn:singular} and the fact that
$\sigma_{\sum_{i=1}^n c_i X_i}^\alpha = \int_E \left\vert
\sum_{i=1}^n c_i f_i(x) \right\vert^\alpha m(dx)$, see relation
(\ref{eq:scale}).

\begin{remark}
A more universal criterion of singularity for stable random
vectors  is in terms of their spectral measure. If  measure
$\Gamma(ds)$ on $\mathbb{S}^{n-1}$ is a spectral measure of an
$\alpha$-stable vector $\vec{X}$ in $\mathbb{R}^n$ and is
concentrated on the intersection of $\mathbb{S}^{n-1}$ with a
$(n-1)$--dimensional linear subspace, then the random vector
$\vec{X}$ is singular. Otherwise, $\vec{X}$  is full--dimensional. For the proof see \cite{KSS13}.
\end{remark}


\subsubsection{Stable Random Fields with an Integral Spectral Representation}

\begin{definition} A random field  $X$  is called {\it{$\alpha$-stable}}\index{stable random field} if all its finite-dimensional distributions
are  $\alpha$-stable.
\end{definition}
Consider  random fields $X=\{ X(t),\, t\in\R^d \}$ of the form
\begin{equation}
 X(t) = \int_E f_t(x) \, M(dx), \quad t \in \mathbb{R}^d, \label{eq:integral_representation}
\end{equation}
where $f_t:E \to \mathbb{R}$ are measurable functions such that $\int_E \vert f_t(x) \vert^\alpha m(dx) < \infty$
and in the case $\alpha=1$ additionally $\int_E \vert f(x) \beta(x) \ln \vert f(x)\vert \vert m(dx) < \infty$ for any $t\in\R^d$. Here
$M$ is an $\alpha$-stable random measure on $(E, \mathcal{E}_0 )$ with control measure $m$ and skewness function $\beta$.
Obviously, the marginals of the random field $X$ in (\ref{eq:integral_representation}) are $\alpha$-stable.
If $\beta(x)=0$ for all $x \in E$ then all finite-dimensional distributions of $X$ are symmetric $\alpha$-stable, so we call $X$ to be a
$S\alpha S$ random field.

A natural question is which stable fields allow for an integral representation (\ref{eq:integral_representation}).
A necessary and sufficient condition for this is the condition of separability of $X$ in probability, see  \cite[Theorem 13.2.1]{ST94}.

\begin{definition}
A stable random field $X=\{X(t),\ t\in M\}$, $M\subseteq\mathbb{R}^d$ is \emph{separable in probability}\index{separable in probability random field} if
there exists a countable subset $M_0\subseteq M$ such that  for every $t\in M$ and any sequence $\{t_k\}_{k\in\N}\subset M_0$ with
$t_k\to t$ as $k\to\infty$ it holds $X(t)=\mbox{plim}_{k\to\infty} X(t_k)$.
\end{definition}
In particular, all stochastically continuous $\alpha$-stable random fields are separable in probability.


\subsection{Dependence Measures for Stable Random Fields}
\label{subsec:2.4}

The dependence of two $\alpha$-stable random variables cannot be digitized by using the covariance because of the absence of the second moments if $\alpha<2$.
We consider two different ways of measuring the degree of dependence of two stable random variables.


\paragraph{\bf Covariation}

\begin{definition}\label{def:covariation1}
Let  $\boldsymbol{X}=(X_1,X_2)^\top$ be an  $\alpha$-stable random vector with $\alpha \in (1,2]$ and spectral measure  $\Gamma$.
 The \emph{covariation}\index{covariation} of $X_1$ on $X_2$ is the real number
$$\left[X_1,X_2\right]_\alpha =  \int_{\mathbb{S}^1} s_1 s_2^{<\alpha-1>} \Gamma(ds_1,ds_2).$$
\end{definition}
It has the following properties.
\begin{theorem}[Properties of Covariation]
Let  $(X_1,X_2, X_3)^\top$ be an $\alpha$-stable random vector with $\alpha \in (1,2]$.
\begin{enumerate}
\item \emph{Linearity in the first entry}:
for $a, b \in \mathbb{R}$ it holds
$$[aX_1+bX_2, X_3]_{\alpha}=a[X_1,X_3]_{\alpha}+b[X_2,X_3]_{\alpha}.$$
\item If $X_1$ and $X_2$ are independent  then $[X_1, X_2]_{\alpha}=0$.
\item \emph{Gaussian case}: for $\alpha=2$, it holds
 $\left[X_1,X_2\right]_2 = 1/2 \cdot \cov(X_1,X_2)$.
\item \emph{Covariation and mixed moments}:
Let $1<\alpha<2$ and $\Gamma$ be spectral measure of  $(X_1,X_2)^\top$ with  $X_1 \sim S_\alpha(\sigma_1, \beta_1,0)$,
$X_2 \sim S_\alpha(\sigma_2,\beta_2,0)$. For $1\leq p<\alpha$, it holds
\begin{equation}
\frac{\mathbb{E}\left(X_1 X_2^{<p-1>}\right)}{\mathbb{E} |X_2|^p}=\frac{[X_1,X_2]_{\alpha}(1-c\cdot \beta_2)+ c\cdot (X_1,X_2)_{\alpha}}{\sigma_2^{\alpha}}, \label{eq:mixed_moments}
\end{equation}
where $(X_1,X_2)_\alpha= \int_{\mathbb{S}^1} s_1\vert s_2\vert^{\alpha-1}\Gamma(ds)$ and
\begin{equation*}
c= \frac{\tan (\alpha \pi/2)}{ 1+\beta_2^2  \tan^2 (\alpha \pi/2 )} \left[ \beta_2 \tan(\alpha \pi/2)-\tan \left(\frac{p}{\alpha} \arctan(\beta_2  \tan (\alpha \pi/2))\right)\right].
\end{equation*}
\end{enumerate}
\end{theorem}
\begin{proof}
 \begin{enumerate}
\item Linearity in the first argument is obvious. However, the covariation is not symmetric, so that there is no linearity in the second argument.
\item To see this, use Exercise~\ref{exe:srv}.
\item The assertion (together with the useful relation  $\var X_i=2\int_{\mathbb{S}^1} s_i^2 \Gamma(d s_1, d s_2)$, $i=1,2$)
follows from the comparison of the characteristic function $\varphi(\theta)$ of the Gaussian random vector  $(X_1, X_2)^\top$  in representations
\eqref{eq:gaus_vect} and \eqref{eq:gaus_vect1}.
\item  See \cite{KSS12} for the proof.
\end{enumerate}
\end{proof}
\begin{remark}
If $X_2$ is symmetric, i.e. $\beta_2=0$, then $c=0$ and formula \eqref{eq:mixed_moments} has the following simple form
\begin{equation*}
\frac{\mathbb{E}\left(X_1 X_2^{<p-1>}\right)}{\mathbb{E} |X_2|^p}=\frac{[X_1,X_2]_{\alpha}}{\sigma_2^{\alpha}},
\end{equation*}
which allows for the estimation of $[X_1,X_2]_{\alpha}$ via empirical mixed moments of $X_1$ and $X_2$.
\end{remark}

For a stable random field $X$ with integral representation (\ref{eq:integral_representation}),
the covariation  can be calculated by the formula
\begin{equation}\label{eq:cov_integral}
[X(t_1), X(t_2)]_{\alpha}=\int_E f_{t_1}(x) (f_{t_2}(x))^{<\alpha-1>} dm(x).
\end{equation}
Notice that its proof given in \cite[Proposition 3.5.2]{ST94} for the $S\alpha S$ case holds true for skewed random fields as well.

\paragraph{\bf Codifference}

Drawbacks of the covariation are the lack of symmetry and the
impossibility to define it for  $\alpha \in (0,1]$. The following
measure of dependence does not have these drawbacks.  That is
however compensated by a mathematically less convenient form.

\begin{definition}\label{dfn:codiff}
Let $(X_1,X_2)^\top$ be an $\alpha$-stable vector. The
\textit{codifference}\index{codifference} of $X_1$ and $X_2$ is
$$\tau (X_1,X_2)= \sigma_{X_1}+\sigma_{X_2}- \sigma_{X_1-X_2},$$
 where $\sigma_{Y}$ is the scale parameter of a $\alpha$-stable random variable $Y$.
\end{definition}

\begin{theorem}[Properties of Codifference]\label{lemma:prop_cd}
\begin{enumerate}
\item \emph{Symmetry}: $\tau(X_1,X_2)=\tau(X_2,X_1)$.
\item If $X_1$ and $X_2$ are independent then $\tau(X_1,X_2)=0$. The inverse statement holds only for  $\alpha\in (0,1)$.
\item \emph{Gaussian case}: for $\alpha=2$, it holds $\tau(X_1,X_2)=\cov(X_1,X_2)$.

\item
Let $(X_1,X_2)$ and $(X'_1,X'_2)$ be S$\alpha$S vectors such that $\sigma_{X_1}=\sigma_{X_2}=\sigma_{X'_1}=\sigma_{X'_2}$.
If $\tau(X_1,X_2)\leq \tau(X'_1,X'_2)$ then for any $c>0$
$$
P\{|X_1-X_2|>c\}\geq P\{|X'_1-X'_2|>c\},
$$
i.e., the larger the codifference, the greater the dependence.
\end{enumerate}
\end{theorem}

\begin{proof}
\begin{enumerate}
\item Symmetry is obvious.
\item Use Exercise~\ref{exe: ind} to see the first part of the statement.
Now let $\tau(X_1,X_2)=0$. It holds $\sigma_{X_1}+\sigma_{X_2}= \sigma_{X_1-X_2}$ iff
$$\int_{\mathbb{S}^1} |s_1|^{\alpha} \Gamma(ds)+\int_{\mathbb{S}^1} |s_2|^{\alpha} \Gamma(ds)=\int_{\mathbb{S}^1} |s_1-s_2|^{\alpha} \Gamma(ds).$$
 We know however that  $|s_1-s_2|^{\alpha}=|s_1|^{\alpha}+|s_2|^{\alpha}$ iff $\alpha<1$ and $s_1 s_2=0$.
\item  It holds $\tau(X_1,X_2)= 1/2( \var \, X_1+ \var \,  X_2- \var (X_1-X_2))= \cov(X_1,X_2)$.
\item See \cite[Property 2.10.6]{ST94} for the proof.
\end{enumerate}
\end{proof}




\subsection{Examples of Stable Processes and Fields}
\label{subsec:2.5}

\paragraph{\bf 1. Stable L\'evy Process} \index{stable L\'evy process}
This is a process defined by $X(t)= M\left([0,t]\right)$,
$t\in\R_+$ where $M$ is an $\alpha$-stable measure on $\R_+$ with
skewness function $\beta$ and Lebesgue control measure multiplied
by $\sigma>0$. $X$ has representation
\eqref{eq:integral_representation} with $f_t(x)=\ind ( x\in
[0,t])$. It obviously holds $X(0)=0$ a.s. Moreover, $X$ has
independent and stationary increments.

Depending on $\beta$ the skewness of the process may vary. So, for
$\alpha<1$ and $\beta\equiv 1$ we obtain a stable L\'evy process
with non-decreasing sample paths, the so--called \emph{stable
subordinator}\index{stable subordinator}. To see this use
one-to-one correspondence between the infinitely divisible
distributions and the  L\'evy processes, thus $X(1)\sim
S_{\alpha}(\sigma, 1,0)$ corresponds to a L\'evy process with the
triplet $(0,0, \frac{\sigma\alpha}{\Gamma(1-\alpha)
cos(\pi\alpha/2)} \frac{dx}{x^{\alpha+1}}\ind (x>0))$, which has
only positive integrable jumps, see also \cite[Examples~21.7 and
24.12]{Sat99}.

\begin{figure}[ht!]

\subfigure[\!\!\! Bisquare \! kernel \! \mbox{$f(x)=\frac{15}{16}
\left(1-\|x\|^2 \right)^2 \ind \left(x\in B_1(0)\right)$.}]{
\label{fig:5}
\includegraphics[scale=.45]{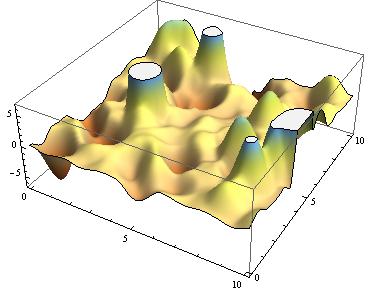}
} \hfill \subfigure[Cylindric kernel $f(x)=\ind \left(x\in
B_1(0)\right)$, $x\in\R^2$.]{
\includegraphics[scale=.45]{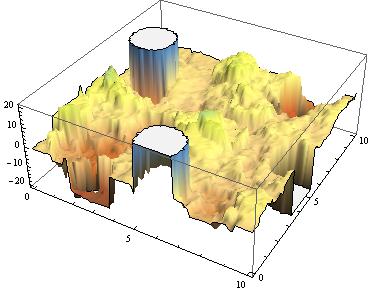}
\label{fig:6}
}
\caption{Continuous (left) and discontinuous (right) realization of a $0.8$-stable moving average random field with $S\alpha S$ random measure $M$.}
%
%
\end{figure}

\paragraph{\bf 2. Stable Moving Average Random Fields} \index{stable moving average random field}
A \emph{stable moving average} $X=\{ X(t), \; t\in\R^d  \}$ is
defined by the formula
$$
X(t)=\int_{\R^d} f(t-s) M(ds), \quad t\in \R^d,
$$
 where  $f\in L^{\alpha}(\R^d)$ is called a \emph{kernel function}\index{kernel function} and $M$ is an $\alpha$--stable random  measure with Lebesgue control measure.
It can be easily seen that $X$ is strictly stationary. See Figures~\ref{fig:5} and \ref{fig:6} for simulated realizations of moving averages in $d=2$ with
the  bisquare and the cylindric kernels.

\emph{Stable Ornstein--Uhlenbeck process} \index{stable
Ornstein--Uhlenbeck process} is a stable moving average process
$X(t)=\int_{-\infty}^t e^{-\lambda(t-s)} M(ds)$, $t\in \R$ where
$M$ is a $S\alpha S$ random measure on $\R$ with the Lebesgue
control measure. The process $X=\{ X(t), \; t\in\R\}$ is strictly
stationary.

\paragraph{\bf 3. Linear Multifractional Stable Motion} \index{linear multifractional stable motion}
is given by $$X(t)=\int_{\R} ((t-x)_+^{H(t)-1/\alpha}-(-x)_+^{H(t)-1/\alpha})M(dx), \quad t\in\R,$$
 where $M$ is an $\alpha$-stable random measure with skewness function $\beta$ and Lebesgue control measure, $\alpha\in (0,2]$.
 The continuous function $H: \R^d \to  (0,1)$ is called a \emph{local scaling exponent}, and $(x)_+=max\{x,0\}$.
It is known that $X$ is a locally self--similar random field, for
more details see  e.g. \cite{ST04,ST04_1}.
 In case $\alpha=2$  we have a Gaussian process called \emph{multifractional Brownian motion}, cf. \cite{PL95}.
For constant $H\in(0,1)$, we get the usual \emph{linear fractional
stable motion} which has stationary increments
 and is $H$--self--similar (see \cite[Sect.~9.5]{BS12}).

\paragraph{\bf 4. Stable Riemann--Liouville Process} \index{stable Riemann--Liouville process}
It is given by $R^H(t)=\int_0^t (t-s)^{H-1/\alpha} M(ds)$, $t\in
\R_+$, where $M$ is an $\alpha$-stable random measure on $\R_+$
 and $H>0$. This is a family of $H$--self--similar random processes. Notice that $R^H$ has no stationary increments, unless $H= 1/\alpha$.
For $\alpha=2$ we get the Gaussian Riemann--Liouville process, see
e.g. \cite[Example~3.4]{Lif12}.

\paragraph{\bf 5. Sub--Gaussian Random Fields} \index{sub--Gaussian random field}
are fields $X$ of the form
$$X\stackrel{d}{=}\{A^{1/2}G(t),\; t \in \R^d\},$$
where
$A \sim S_{\alpha/2}((\cos(\pi\alpha/4))^{2/\alpha},1,0)$ and $G=\{G(t), \; t \in \R^d\}$ is a
a zero mean Gaussian random field with a positive definite covariance function which is independent of $A$. The following
lemma (cf. \cite[Proposition~3.8.1]{ST94}) shows that $X$ is $\alpha$-stable.
\begin{lemma} If the random variable $A$ is as above
and  $\xi\sim N(0,2\sigma^2)$ independent of $A$ then $X=A^{1/2}\xi\sim S_{\alpha}(\sigma,0,0)$.
\end{lemma}
To prove the lemma, it suffices to calculate the characteristic
function of $X$ using the conditional expectation provided that
$A$ is fixed. If $G$ is stationary then the resulting
sub--Gaussian field $X$ is strictly stationary as well.

A strictly stationary sub--Gaussian random field $X$ with a mean
square continuous Gaussian component $G$ is not ergodic since it
differs from $G$ by a random scaling. A sufficient condition for
ergodicity of $G$ is that its spectral measure has no atoms, see
\cite[Theorem A]{Tem73}.


\section{Extrapolation of Stationary Random Fields}
\label{sec:3}

Let $X=\{  X(t),\; t\in\R^d\}$ be a stationary (in the appropriate
sense to be specified later) random field. We are looking for a
\emph{linear predictor}\index{linear predictor} $\widehat{X}(t)$
of the unknown field value $X(t)$ at location $t\in \R^d$ based on
observations $X(t_1),\ldots,X(t_n)$ at locations $t_1,\ldots,t_n$,
$n \in \N$ of the form
\begin{equation}
\widehat{X}(t) = \sum_{i=1}^n \lambda_i (t) X(t_i) +\lambda_0 (t). \label{eq:linear_estimator2}
\end{equation}
The weights $\lambda_0(\cdot),\ldots,\lambda_n(\cdot)$ are functions of $t,t_1,\ldots,t_n$ which may depend on the distribution of $X$.
For simplicity of notation, we omit all their arguments except for $t$. They have to be computed in a way (which
depends on the integrability properties of $X$) such that the predictor is in some regard close to $X(t)$.

\begin{definition}
A predictor $\widehat{X}(t)$ for $X(t)$ is called
\begin{enumerate}
\item \emph{exact}\index{predictor!exact}  if $\widehat{X}(t) = X(t)$ a.s. whenever $t = t_i$ for any $i \in \{1,\ldots,n\}$. In this case, the predictor
 $\hat{X}(\cdot)$ is an extrapolation surface for $X(\cdot)$ with knots $t_1, \ldots, t_n$.

\item  \emph{unbiased}\index{predictor!unbiased} if $\E |X(0)|<\infty$ and $\mathbb{E} (\widehat{X}(t) - X(t)) = 0$.

 \item \emph{continuous}\index{predictor!continuous} if weights $\lambda_i(\cdot)$, $i=0,\ldots,n$ are continuous with respect to $t$, i.e.,
any realization of  $\hat{X}$ is continuous in $t\in\R^d$.
\end{enumerate}
\end{definition}

\subsection{Kriging Methods for Square Integrable Random Fields}

If the field $X$ has finite second moments then the most popular
prediction technique for $X$ in geostatistics is the so--called
\emph{kriging}\index{kriging}. It is named after D.G. Krige who
first applied it (in 1951) to gold mining. Namely, he predicted
the  size of a gold deposit by collecting the data of gold
concentration at some isolated locations. Apart from kriging,
there are many other prediction techniques such as  inverse
distance, spline and nearest neighbor interpolation,
triangulation, see for details \cite[Sect.~5.9.2]{Cres91},
\cite[Chapt.~3]{WO07}, \cite{Sib81}, etc. However, the latter
methods ignore the correlation structure contained in the spatial
data; see \cite[Sect.~3.4.5, p.180; Chapt. 5.9]{Cres91},
\cite{Fra82}, \cite{Las94} for their comparison.

The main idea of  kriging is to compute prediction weights $\lambda_i$ by minimizing
 the mean square error between the predictor and the field itself, i.e., solve the minimization problem
\begin{equation}\label{eq: LSM}
\E (X(t)-\hat X(t))^2\rightarrow \min_{\lambda_0, \ldots\lambda_n}
\end{equation}
under some additional conditions on $\lambda_i$ for each fixed $t\in \R^d$.

Depending on the assumptions about $X$,  numerous variants of kriging are avaliable. We mention just few of them and
refer an interested reader to the vast literature.
\begin{enumerate}
\item \emph{Simple kriging}: for square integrable random  fields
$X$ with known mean function $\E X(t)=m(t)$, $t\in\R^d$. See
Section~\ref{subsec:3.1}.

\item  \emph{Ordinary kriging}: for second order intrinsic
stationary  random fields  $X$ (with unknown but constant mean).
See Section~\ref{subsec:3.2}.

\item \emph{Kriging with drift}:  $\E X(t)=a+b \|t\|$, $a,\ b\in
\R$ and these constants are unknown. See \cite[Sect. 3.4.6]{CD99}
for details.

\item \emph{Universal kriging}: the unknown mean $\E X(t)=m(t) \neq const$ belongs to some parametric family of functions, see \cite{CD99,Wack03}.
 Ordinary kriging and kriging with drift are special cases of universal  kriging.
\end{enumerate}

\subsection{\index{simple kriging}{Simple Kriging}}
\label{subsec:3.1} Let $X$ be a square integrable random field
with known mean function $m(t)$. It is easy to see that the
minimum of the mean square error
\begin{equation*}
\E (X(t)-\hat X(t))^2= \var (X(t)-\hat X(t))+(\E(X(t)-\hat X(t)))^2
\end{equation*}
is attained exactly  when the predictor $\hat X(t)$ is unbiased, i.e. if $\E \hat{X}(t)=\E X(t)$. This yields $\lambda_0(t)=m(t)- \sum_{i=1}^n \lambda_i(t)m(t_i)$ and
\begin{equation*}
\widehat{X}(t) = \sum_{i=1}^n \lambda_i (t) (X(t_i) - m(t_i))+m(t).
\end{equation*}
It follows from the above relation that the knowledge of function
$m$ leads to centering the field $X$ (subtracting $m$) in the
prediction.

 Taking derivatives of the goal function in (\ref{eq: LSM}) with
respect to $\lambda_i$, we obtain
\begin{equation}\label{eq:system_kriging}
\sum_{i=1}^n \lambda_i(t) \cov(X(t_i), X(t_j))= \cov(X(t), X(t_j)),\   j=1,\ldots,n.
\end{equation}
The matrix form of this system of equations is
$$
\Sigma\cdot \lambda(t)=\sigma(t),
$$
where $\Sigma=[\cov(X(t_i), X(t_j))]_{i,j=1}^n$ is the covariance
matrix, $\lambda(t)=(\lambda_1(t), \ldots, \lambda_n(t))^\top$,
$$
\sigma(t)= (\cov(X(t), X(t_1)), \ldots , \cov(X(t), X(t_n)))^\top.
$$
If $\Sigma$ is non-degenerate then the solution exists and is
unique. The covariance matrix is non-degenerate if the covariance
function of $X$ is positive definite and all $t_i$, $i=1,\ldots,n$
are distinct.

\begin{exercise}
Let the random field  $X=\{X(t),\  t\in\R^d\}$ be as above. Show
that the random vector $(X(t_1), \ldots, X(t_n))^\top$ is singular
iff $\det \Sigma=0$. {\it{Hint:}} A symmetric matrix is positive
definite (positive semi--definite) if and only if all of its
eigenvalues are positive (non--negative).
\end{exercise}

Finally, we have the following form of the predictor:
\begin{equation}\label{eq:matrix_krig}
\hat{X}(t)=\bar{X}^\top\Sigma^{-1}\sigma(t),
\end{equation}
where $\bar{X}=(X(t_1), \ldots, X(t_n))^\top$.

Let $\delta_{ij} = \ind (i=j)$ be the Kronecker delta.

\paragraph{\bf Properties of Simple Kriging}

\begin{enumerate}
\item {\bf Exactness}: to see that $\hat{X}(t_j)=X(t_j)$ for any
$j$, set $t=t_j$  and check that $\lambda_i(t_j)=\delta_{ij}$,
$i,j=1,\ldots, n$ is the solution of system of equations
(\ref{eq:system_kriging}). \item {\bf Continuity and smoothness}:
rewrite (\ref{eq:matrix_krig}) as $ \hat{X}(t)=b^\top \sigma(t)$
with $b=\Sigma^{-1}\bar{X}$ which means that sample path
properties of the extrapolation surface such as continuity and
smoothness directly depend on the properties of $\sigma(t)$.
Thus if the covariance function is continuous and smooth, so is the extrapolation surface. See Figure~\ref{fig:11}.\\
\item {\bf Shrinkage property}: The mean prediction error $\E
(\hat{X}(t)-X(t))^2$ can be found by direct calculations using the
system (\ref{eq:system_kriging}). Thus
\begin{equation}\label{eq:smoothing}
\E (\hat{X}(t)-X(t))^2=\var X(t)-\var \hat{X}(t).
\end{equation}
Equation (\ref{eq:smoothing}) yields the following
\index{shrinkage property}{shrinkage property}: for all $t\in
\R^d$
\begin{equation}\label{eq:shrinkKr}
\var \hat{X}(t)\leq \var X(t).
\end{equation}
The simple kriging predictor is less dispersed than the data. In a
sense, kriging performs linear averaging (or smoothing) and does
not perfectly imitate the trajectory properties of the original
random field.

\item {\bf Geometric interpretation}: The predictor $\hat{X}(t)$
for any fixed $t$ can be interpreted as a metric projection of
$X(t)$ onto the linear subspace $L_n=\mbox{span} \{X(t_1), \ldots,
X(t_n)\}$ of Hilbert space $L^2(\Omega, {\cal F}, \P)$ with scalar
product $\langle X,Y\rangle=\E (X Y)$ for $X,Y\in L^2(\Omega,
{\cal F}, \P)$, that is,
\begin{equation}
\hat{X}(t)=\mbox{Proj}_{L_n} X(t)=\argmin_{\xi\in L_n} \langle
X(t)-\xi, X(t)-\xi \rangle.
\end{equation}
It is known from the Hilbert space theory that this projection is
unique if the vector $(X(t_1), \ldots, X(t_n))^\top$ is not
singular (cf. Definition~\ref{dfn:singular}).

\item {\bf Orthogonality}: The above projection is also
orthogonal, i.e., $\langle \hat{X}(t)-X(t), \xi \rangle=0$ for all
$\xi \in L_n$. In particular, it holds
\begin{equation}\label{eq:orthogonality}
\langle \hat{X}(t)-X(t), X(t_i) \rangle=0\ \ \text{for all}\ \
i=1,\ldots, n
\end{equation}
which rewrites as a dependence relation
\begin{equation*}
\E \left( \hat{X}(t) X(t_i)\right) =\E \left( X(t) X(t_i) \right)\
\ \text{for all}\ \ i=1,\ldots, n
\end{equation*}
yielding
\begin{equation*}
\cov(\hat{X}(t)-X(t), \hat{X}(s))=0, \quad s\in \mathbb{R}^d.
\end{equation*}
\begin{exercise}
Prove  relation \eqref{eq:smoothing} via the Pythagorean theorem.
\end{exercise}

\item {\bf Gaussian case}: Under the assumptions that $X$ is
Gaussian and $\Sigma$ is non--singular it is easy to show that
\begin{equation}\label{eq:GaussCondExpSimpleKrig}
\hat{X}(t)= \E \left(X(t)| X(t_1), \ldots, X(t_n)\right),\quad
t\in\R^d.\end{equation}
\begin{exercise}Prove relation \eqref{eq:GaussCondExpSimpleKrig}
 using the uniqueness of the kriging predictor and
the following properties of the conditional expectation
and of the Gaussian multivariate distribution, respectively:
\begin{enumerate}
\item[1)] $\E ((\eta-\E(\eta|\xi))h(\xi))=0$ for random variables
$\xi,\eta$ and any measurable function $h(\cdot)$ , \item[2)] If
$\eta, \xi_1,\ldots, \xi_n$ are jointly Gaussian then there exist
real numbers $\{a_i\}_{i=1}^n$ such that $\E(\eta| \xi_1, \ldots,
\xi_n)=\sum_{i=1}^n a_i \xi_i$. \end{enumerate}
\end{exercise}

In the Gaussian case, simple kriging has additional properties of
\begin{enumerate}
\item {\bf Conditional unbiasedness}:\index{simple
kriging!conditional unbiasedness}
$\E\big(X(t)|\hat{X}(t)\big)=\hat{X}(t)$ a.s. for any $t\in\R^d$,
cf. \cite[p.~164]{CD99}. This property is important in practice
for resource assessment problems and selective mining.
 \item {\bf Homoscedasticity}\index{simple
kriging!homoscedasticity}: The conditional mean square estimation
error does not depend on the data, i.e.,
$$
\E\left(\big(\hat{X}(t)-X(t)\big)^2| X(t_1), \ldots,
X(t_n)\right)= \E \big(\hat{X}(t)-X(t)\big)^2 \ \mbox{ a.s. for
any } \ t\in\R^d.
$$
\end{enumerate}
\end{enumerate}

\subsection{Ordinary Kriging}
\label{subsec:3.2}

When the mean $m$ of a square integrable random field $X$ is
constant but unknown the \emph{ordinary kriging}\index{ordinary
kriging} can be applied. We are looking for a predictor in the
form (\ref{eq:linear_estimator2}). For an arbitrary (but fixed)
location $t\in\R^d$, the mean square prediction error is
$$
\E
\left(\hat{X}(t)-X(t)\right)^2=\var(\hat{X}(t)-X(t))+\left(\lambda_0+\left(\sum_{i=1}^n
\lambda_i-1\right) m\right)^2.
$$
 Assuming that
\begin{equation}\label{eq:OK_constr}
 \lambda_0=0, \quad \sum_{i=1}^n \lambda_i =1,
\end{equation}
we get the smallest possible error together with unbiasedness $\E
\hat{X}(t)=\E X(t)$. The ordinary kriging predictor writes then
$$
\hat{X}(t)=\sum_{i=1}^n \lambda_i X(t_i),\quad t\in\R^d.
$$

The prediction error can be computed as
$$
\E \left(\hat{X}(t)-X(t)\right)^2=\sum_{i,j=1}^n
\lambda_i\lambda_j \cov(X(t_i), X(t_j))-2\sum_{i=1}^n \lambda_i
\cov(X(t_i), X(t))+\var X(t).
$$
One should minimize this error under the constraint
(\ref{eq:OK_constr}).

Taking partial derivatives of the  Lagrange function
$$
L(\vec{\lambda}, \mu)=\E
(\hat{X}(t)-X(t))^2+2\mu\left(\sum_{i=1}^n \lambda_i -1\right)
$$
with respect to $\lambda_i=\lambda_i(t)$, $i=1,\ldots, n$, and
$\mu=\mu(t)$ and putting them equal to zero we obtain the
following system of $n+1$ linear equations
\begin{eqnarray*}
\begin{cases}
 \sum_{i=1}^n\lambda_i\cov(X(t_i), X(t_j))+\mu= \cov(X(t_j), X(t)), \quad j=1,\ldots,n,\\
 \sum_{i=1}^n \lambda_i=1
\end{cases}
\end{eqnarray*}
for each $t\in\R^d$ of interest. The solution
$(\lambda_1,\ldots,\lambda_n,\mu)^\top$ of this system is unique
iff the covariance matrix of the vector $\big(X(t_1),\ldots,
X(t_n)\big)^\top$ is non--singular.

The above linear system of equations can be rewritten in terms of
 variogram $\gamma(\cdot, \cdot)$. By formula
(\ref{eq:variogram}) and direct calculation we get the following
ordinary kriging system of equations with respect to the weights
$\lambda_i$, $i=1,\ldots, n$, and $\mu$:
\begin{eqnarray*}
\begin{cases}\sum\limits_{i=1}^n \lambda_i \gamma(t_i, t_j)+\mu=\gamma(t_j, t), \quad j=1,\ldots,n,\\
\sum\limits_{i=1}^n \lambda_i=1.
\end{cases}
\end{eqnarray*}
The corresponding mean square prediction error is
$$
\sigma^2_{OK}=\E \big(\hat{X}(t)-X(t)\big)^2=\sum_{i=1}^n
\lambda_i \gamma(t_i, t)+\mu.
$$
\begin{exercise}
Show that $\mu=-(1-{\bf e}^\top \Gamma^{-1}\gamma)/{\bf e}^\top
\Gamma^{-1}{\bf e}$, where ${\bf e}$ is the unit vector,
$\gamma=(\gamma(t_1, t), \ldots, \gamma(t_n, t))^\top$   and
$\Gamma=[\gamma(t_i, t_j)]_{i,j=1,\ldots, n}$.
\end{exercise}
The main advantage of this way of posing the problem is that it is
solvable even if the variance of $X(t)$ is infinite whereas the
variogram is finite, e.g., if $X$ is intrinsic stationary of order
two.

\paragraph{\bf Properties of the Ordinary Kriging}
\begin{enumerate}
\item {\bf Exactness}: For $t=t_j$, notice that
$\lambda_i(t_j)=\delta_{ij}$, $i,j=1,\ldots, n$, $\mu(t_j)=0$ is a
solution of the ordinary kriging system.

\item {\bf Orthogonality}: For any real weights $a_i$,
$i=1,\ldots,n$ with the property $\sum_{i=1}^n a_i=1$ it holds
$$
\left\langle \hat{X}(t)-X(t), \sum_{i=1}^n a_i
X(t_i)\right\rangle=0.
$$
\item {\bf Conditional unbiasedness}: The ordinary kriging
predictor  reduces the conditional bias
$\E\left(X(t)|\hat{X}(t)\right) -\hat{X}(t)$. To see this, check
the following formula showing that the minimum of the kriging
error corresponds to the minimum of the conditional bias error:
\begin{equation*}
\E\left(\E\left(X(t)|\hat{X}(t)\right) -\hat{X}(t)\right)^2= \E
\left(\hat{X}(t)-X(t)\right)^2-\E
\left(\var\left(X(t)|\hat{X}(t)\right)\right),
\end{equation*}
cf. \cite[p.185]{CD99}. For the proof of this formula, the
following \emph{ law of total variance}  is used
\begin{equation*}
\var Y=\var \left(\E(Y|Z)\right)+\E\left( \var(Y|Z)\right)
\end{equation*}
 as well as $\E \left(Y\cdot \E(Z|Y) \right)=\E (YZ)$ for any random variables $Y,$ $Z$ defined on the same probability space.
\end{enumerate}

\begin{example}\label{ex:OrdKrig}
A simulated  realization (see Figure~\ref{fig:9}) of a centered
stationary isotropic Gaussian random field $X=\{ X(t), \;
t\in[0,10]^2 \}$ with Whittle--Mat\'ern--type covariance function
$C(s,t)=2\ind (s=t)+2\ind (s\neq t) \| s-t\| K_1(2 \| s-t\|)$
exhibiting a nugget effect of height one is observed on a grid of
locations $\{ (3i,2j), \; i,j\in \N \cap [0,3] \}$. The
corresponding theoretical variogram together with the Matheron
estimator (given in \cite[Formula (9.67)]{BS12}) are shown on
Figure~\ref{fig:10}. A Whittle--Mat\'ern--type variogram model
with a nugget effect $\sigma ^2$
$$
\gamma(s,t)=\ind(s\neq t) \left( \sigma^2+b -b 2^{1-\nu} (a\|s-t
\|)^{\nu}K_{\nu}(a \|s-t \|) \right), \quad s,t\in\R^d,
$$
was fitted to the estimated variogram by an ordinary least squares
method yielding the parameter estimates $\hat{\sigma}^2 =
0.9327665235438869$, $\hat{a}=1.9674556902269302$,
$\hat{b}=1.0672476194785714$. An extrapolation by ordinary kriging
with the fitted variogram model $\gamma$ is shown on
Figure~\ref{fig:11}.
\begin{figure}[ht!]
\subfigure[Simulated realization of a stationary Gaussian random
field with nugget effect.]{ \label{fig:9}
\includegraphics[scale=.35]{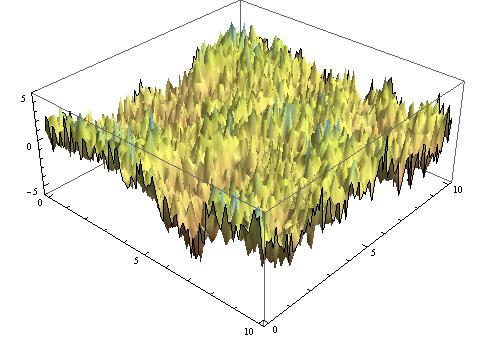}
} \hfill \subfigure[Extrapolation by ordinary kriging for the
field in Figure \ref{fig:9}]{ \label{fig:11}
\includegraphics[scale=.35]{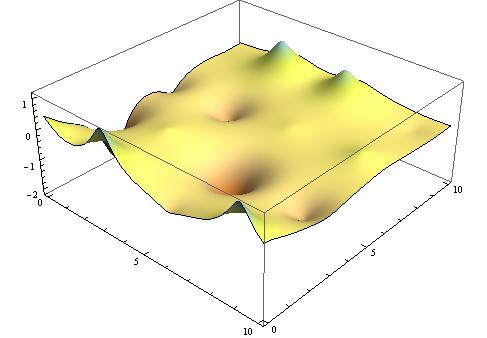}
} \caption{Application of ordinary kriging to simulated data from
Example \ref{ex:OrdKrig}}
\end{figure}
\begin{figure}[ht!]
\sidecaption
\includegraphics[scale=.30]{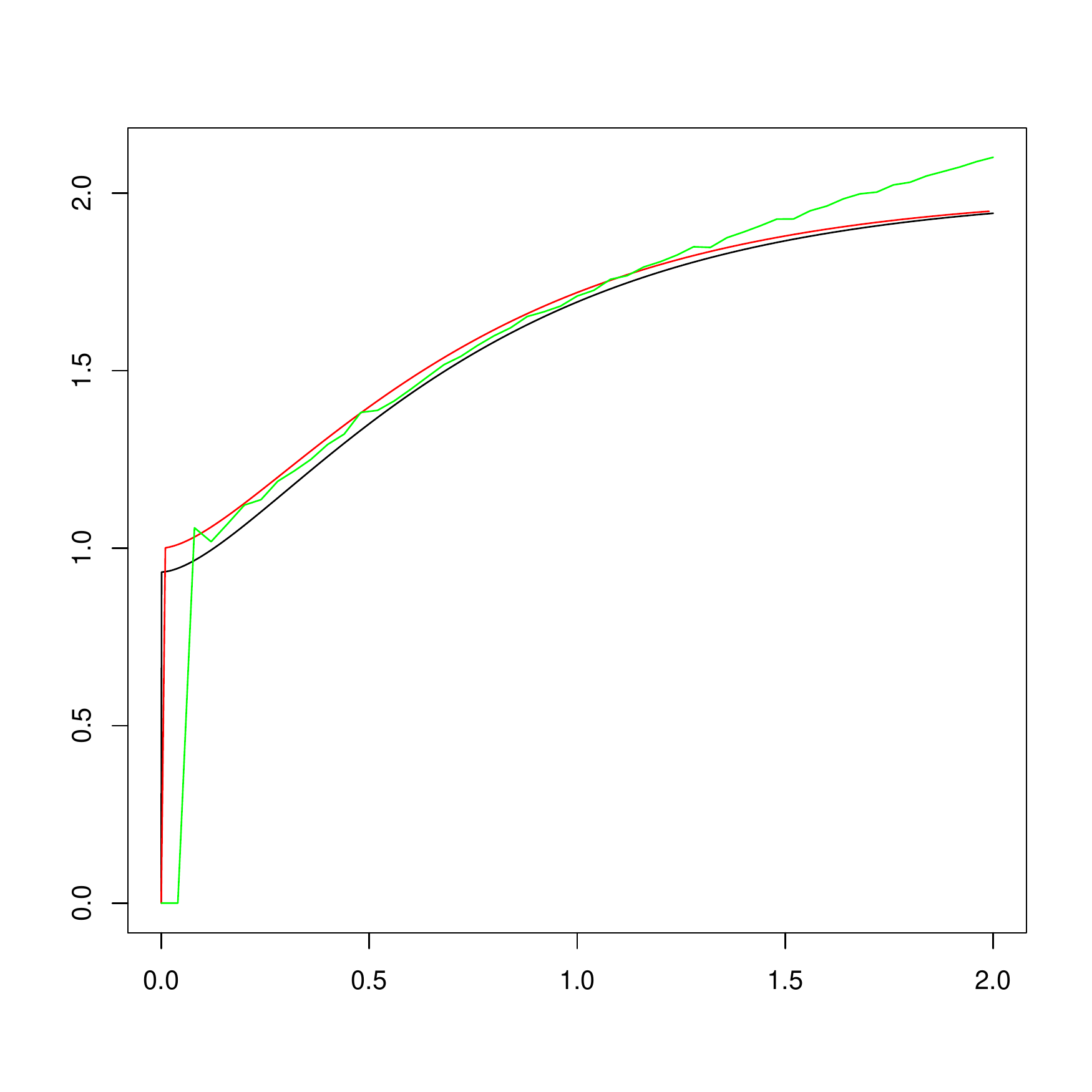}
\caption{Theoretical variogram (red),  estimator (green) and
fitted variogram (black) for the realization in Figure
\ref{fig:9}, compare Example \ref{ex:OrdKrig}} \label{fig:10}
\end{figure}
\end{example}

\section{Extrapolation of Stable Random Fields} \label{sec:4}

Let $X$ be an $\alpha$-stable random field having integral
representation
\begin{equation}\label{eq:RFX}
 X(t) = \int_E f_t(x) M(dx), \quad t \in \mathbb{R}^d,
\end{equation}
confer formula (\ref{eq:integral_representation}). For $\alpha\in
(1,2]$, assume that  the field $X$ is centered. If
$\alpha\in(0,1]$, the mean value of $X$ does not exist.

We are looking for a predictor $\widehat{X}(t)$ of the value
$X(t)$ at location $t \in \mathbb{R}^d$ based on the random vector
$(X(t_1),\ldots,X(t_n))^\top$ in the form
\begin{equation}\label{eq:LinPred}
\widehat{X}(t) = \sum_{i=1}^n \lambda_i X(t_i).
\end{equation}


Let $T_j=\{ t_{j,1},\ldots, t_{j,n_j}  \}$, $j\in\N$ be a sequence
of observation locations such that $\dist(T_j,t)\to 0$ as
$j\to\infty$ where $\dist(A,B)=\inf \{ \|x-y\|:\; x\in A,\, y\in B
\} $ is the Euclidean distance between two arbitrary sets
$A,B\subset \R^d$. The predictor $\hat{X}_j(t)=\sum_{i=1}^{n_j}
\lambda_i^{(j)} X(t_{j,i})$ is \emph{weakly
consistent}\index{predictor!weakly consistent} if
$\hat{X}_j(t)\plongj X(t)$ for any $t\in\R^d$. It is
\emph{stochastically continuous}\index{predictor!stochastically
continuous} if $\hat{X}_j(s)\plongs \hat{X}_j(t)$ for any $j\in\N$
and $t\in\R^d$.

Let
\begin{equation}\label{eq:L_alpha_norm}
\Vert f \Vert_{\alpha} = \left(\int_E \vert f(x)\vert^\alpha
m(dx)\right)^{1/\alpha} \end{equation} denote the norm of $f\in
L^\alpha(E,m)$, $\alpha\ge 1$.
\begin{theorem}\label{th:Consistency}
Let the $\alpha$--stable random field $X$ in \eqref{eq:RFX} be
stochastically continuous, $\alpha\in (1,2]$. Let the  predictor
$\hat{X}_j$ defined above exist and be unique, exact and
stochastically continuous for any $j\in\N$. Then $\hat{X}_j$ is
weakly consistent.
\end{theorem}
\begin{proof}\footnote{The idea of this proof belongs to Adrian Zimmer.}
Fix an arbitrary $t\in\R^d$. By \cite[Proposition 3.5.1]{ST94}, it
is sufficient to show that $\sigma_{\hat{X}_j(t)- X(t)}\to 0$ as
$j\to\infty$ to prove weak consistency. Let $s_j\in T_j$ be the
point at which $\dist(s_j,t)=\dist(T_j,t)$ for any $j\in\N$. It is
clear that $s_j\to t$ as $j\to\infty$. Since $\hat{X}_j$ is exact
it holds $\hat{X}_j(s_j)=X(s_j)$ for any $j$. Then we have
\begin{multline*}
\sigma_{\hat{X}_j(t)- X(t)}= \| \sum_{i=1}^{n_j} \lambda_i^{(j)}
f_{t_{j,i}}-f_t\|_\alpha  \le \| f_{s_j}-f_t\|_\alpha + \|
\sum_{i=1}^{n_j} \lambda_i^{(j)} f_{t_{j,i}} - f_{{s_j}}\|_\alpha
\to 0
\end{multline*}
as $j\to\infty$ by \cite[Proposition 3.5.1]{ST94}), stochastic
continuity of $X$ and $\hat{X}_j$ as well as  exactness of
$\hat{X}_j$.
\end{proof}


\subsection{Least Scale Predictor}
\label{subsec:4.1}

For $\alpha\in(0,2]$, consider the following optimization problem

\begin{equation}
\sigma_{\widehat{X}(t)-X(t)}^\alpha = \int_E \left|f_{t}(x) -
\sum_{i=1}^n \lambda_i f_{t_i}(x)\right|^\alpha m(dx) \, \to \,
\min_{\lambda_1, \ldots, \lambda_n}. \label{eq:LSL}
\end{equation}
It is clear that the solution of the optimization problem (in case
if it exists and is unique) will be an extrapolation. To see this,
put $t=t_j$ and $\lambda_i(t_j)=\delta_{ij}$, $i,j=1,\ldots, n$.

The predictor $\widehat{X}(t)$ based on a solution of this
minimization problem is called \textit{least scale linear (LSL)
predictor}\index{least scale linear (LSL) predictor}. This method
is reminiscent of the least mean square error property (\ref{eq:
LSM}) of the kriging.

If $\alpha\in (1,2]$ it is easy to see that any solution of the
problem
 (\ref{eq:LSL}) is also a solution of the following system of equations
\begin{equation}
\int_E f_{t_j}(x) \left(f_{t}(x) - \sum_{i=1}^n \lambda_i
f_{t_i}(x)\right)^{<\alpha-1>}\hspace*{-0.5cm}m(dx)=0, \quad j=1,
\ldots, n,\label{eq:LSL_system}
\end{equation}
or equivalently
\begin{equation}
\left[X(t_j),X(t) - \sum_{i=1}^n\lambda_i X(t_i)\right]_\alpha=0,
\quad j=1, \ldots , n, \label{eq:LSL_covariation}
\end{equation}
where $[\cdot, \cdot]_{\alpha}$ is the covariation, see Definition
\ref{def:covariation1}.
\begin{exercise}\label{ex:LSL}
 Show that any solution of the problem (\ref{eq:LSL}) solves also the system of equations (\ref{eq:LSL_system}) or (\ref{eq:LSL_covariation}). Use the
 dominated convergence theorem.
\end{exercise}

Notice that equations in (\ref{eq:LSL}) are nonlinear in
$\lambda_1,\ldots,\lambda_n$ if $\alpha < 2$ because the
covariation is not linear in the second argument (cf.
Section~\ref{subsec:2.4}). Thus, numerical methods have to be
applied  to solve problem \eqref{eq:LSL}.

\paragraph{\bf Properties of LSL predictor}

Assume  $1<\alpha\leq 2$. For the case $0<\alpha\leq 1$ see
Section \ref{subsec:4.4}.

\begin{theorem}\label{th:LSL1}
The LSL predictor exists. If the random vector
$\vec{X}=(X(t_1),\ldots,X(t_n))^\top$ is full--dimensional then
the LSL predictor is unique.
\end{theorem}
\begin{proof}
We are using the properties of the best approximation in
$L^\alpha(E,m)$-spaces for $1< \alpha\leq 2$. Let
$L=span\{f_{t_1}, \ldots,  f_{t_n}\}$. This is  a finite
dimensional space. Denote for simplicity $f=f_t$ and
$E(f)=\inf_{x\in L}\|f-x\|_{\alpha}$. Let us show that this
infimum is attained in $L$.

Consider $\{x_m\}_{m\in\mathbb{N}}$  such that $x_m\in L$ $\forall
m\in\mathbb{N}$    and $\|x_m-f\|_{\alpha}\to E(f)$ as $m\to
\infty$. By the triangle inequality  $\|x_m\|_{\alpha}\leq
\|f\|_{\alpha}+\|f-x_m\|_{\alpha}$, so $\{x_m\}_{m\in\mathbb{N}}$
is a bounded sequence in a finite dimensional subspace. Thus,
there exists a convergent subsequence $\{m_j\}_{j\in \mathbb{N}}$
and $f_0\in L$ such that $\|x_{m_j}-f_0\|_{\alpha}\to 0$ as $j\to
\infty$. Since  $\|f-x_{m_j}\|_{\alpha}\to \|f-f_0\|_{\alpha}$ and
$\|f-x_{m_j}\|_{\alpha}\to E(f)$ as $j\to\infty$, it holds
$E(f)=\|f-f_0\|_{\alpha}$. So $f_0$ is the best approximation.

For the proof of uniqueness, we use the strict convexity property.
If $\alpha>1$ the space $L^{\alpha}(E,m)$ is \emph{strictly
convex}\index{strict convexity of $L^\alpha(E,m)$} (see e.g.
\cite[p. 59]{DeVL93}), i.e. for all $g_1, g_2 \in L^{\alpha}(E,m)$
such that $\|g_1\|_{\alpha}=\|g_2\|_{\alpha}=1$, $g_1\neq g_2$  it
follows $\|\beta g_1+(1-\beta) g_2\|_{\alpha}<1$ for any $\beta\in
(0,1)$.

Take  $y_j=\sum_{i=1}^n \lambda_i^{(j)} f_{t_i}\in L,$ $j=1, 2$
such that $y_1\neq y_2$ and
$\|f-y_1\|_{\alpha}=\|f-y_2\|_{\alpha}=E(f)$.  Thus by strict
convexity we have
\begin{equation*}
E(f)\leq \left\|f-\frac{1}{2} \left(y_1+y_2
\right)\right\|_{\alpha} =\left\|\frac{1}{2}
\left(f-y_1\right)+\frac{1}{2}\left(f-y_2\right)\right\|_{\alpha}<E(f).
\end{equation*}
So we obtain a contradiction, and  $y_1=y_2=f_0$. By the
full--dimensionality of the random vector $\vec{X}$  and by
Lemma~\ref{lemma:great_subsphere} one can easily see that  the set
of the weights $\lambda_i$ in the representation $f_0=\sum_{i=1}^n
\lambda_i f_{t_i}$ is unique.
\end{proof}

\begin{theorem}[\cite{KSS12}]\label{th:LSL2}
Let the stable random field $X$ in \eqref{eq:RFX} be
stochastically continuous. If the random vector
$\vec{X}=\big(X(t_1),\ldots,X(t_n)\big)^\top$ is full--dimensional
then the LSL predictor is continuous.
\end{theorem}


\subsection{Covariation Orthogonal Predictor}
\label{subsec:4.2}

Throughout this Section, assume $\alpha\in (1,2]$.  The linear
predictor \eqref{eq:LinPred} with weights
$\lambda_1,\ldots,\lambda_n$ that are a solution of the system of
equations
\begin{equation}
\left[X(t) - \sum_{i=1}^n\lambda_i X(t_i), X(t_j)\right]_\alpha=0,
\quad j=1, \ldots , n \label{eq:COL}
\end{equation}
is called \textit{Covariation Orthogonal Linear (COL)
predictor}\index{Covariation Orthogonal Linear (COL) predictor}.
If the solution of (\ref{eq:COL}) exists and is unique then it is
an exact predictor, since we can put $\lambda_i(t_j)=\delta_{ij}$,
$i,j=1,\ldots, n$. This extrapolation method is reminiscent of the
generic orthogonality property of simple kriging, cf. relation
(\ref{eq:orthogonality}). It is also \emph{symmetric} (in a sense)
to the LSL predictor, compare the systems
\eqref{eq:LSL_covariation} and
 \eqref{eq:COL}. In contrast to
\eqref{eq:LSL_covariation}, the system \eqref{eq:COL} is linear
which makes the computation of the weights $\lambda_i$ easier.

Introduce the \emph{covariation function}\index{covariation!
function} $\kappa: \mathbb{R}^d\times\mathbb{R}^d\to \mathbb{R}$
of $X$ by
\begin{equation}
\kappa(s,t)=[X(s), X(t)]_{\alpha}.\label{eq:covariation_function}
\end{equation}
Note that this function is not symmetric in its arguments, as
opposed to the covariance function, cf. Definition
\ref{dfn:covariance}.

By additivity of the covariation in the first argument (see
Section~\ref{subsec:2.4}), the system (\ref{eq:COL}) rewrites as
\begin{equation}
 \begin{pmatrix}
 \kappa(t_1, t_1) & \cdots & \kappa(t_n, t_1) \\
 \vdots & \ddots & \vdots \\
 \kappa(t_1, t_n) & \cdots & \kappa(t_n, t_n)\end{pmatrix} \begin{pmatrix}
        \lambda_1 \\
        \vdots \\
        \lambda_n
          \end{pmatrix} = \begin{pmatrix}
                \kappa(t, t_1) \\
                \vdots \\
                \kappa(t, t_n)
\end{pmatrix}. \label{eq:COL_sys}
\end{equation}
 If matrix $K=[\kappa(t_i, t_j)]_{i,j=1,\ldots, n}$ is positive definite
the solution of this system exists and is unique.

For moving average and for sub--Gaussian fields $X$, sufficient
conditions for the positive definiteness of $K$ can be given.

\subsubsection{The COL Predictor for Moving Averages}

Consider a moving average stable random field $X$  with
representation
\begin{equation*}
X(t) = \int_{\mathbb{R}^d} f(t-x)\, M(dx), \quad t \in
\mathbb{R}^d,
\end{equation*}
where $M$ is an $\alpha$--stable random measure with Lebesgue
control measure and $f\in L^{\alpha}(\mathbb{R}^d)$ (see
Section~{\ref{subsec:2.5}} for the definition). By strict
stationarity of $X$, it holds $[X(h),X(0)]_\alpha =
[X(t+h),X(t)]_\alpha$ for all $t,\ h \in \mathbb{R}^d$. With
slight abuse of notation, we write $\kappa(s-t)=
[X(s-t),X(0)]_\alpha=\kappa(s,t)$, $s, t\in \mathbb{R}^d$ and the
system of equations (\ref{eq:COL_sys}) is equivalent to
\begin{equation}
 \begin{pmatrix}
 \kappa(0) & \cdots & \kappa(t_n-t_1) \\
 \vdots & \ddots & \vdots \\
 \kappa(t_n-t_1) & \cdots & \kappa(0)\end{pmatrix} \begin{pmatrix}
        \lambda_1 \\
        \vdots \\
        \lambda_n
          \end{pmatrix} = \begin{pmatrix}
                \kappa(t-t_1) \\
                \vdots \\
                \kappa(t-t_n)
\end{pmatrix}.\label{eq:soc4}
\end{equation}
The next theorem gives a sufficient condition for the existence
and uniqueness of the COL predictor.
\begin{theorem}\label{th:COL_MA_pd}
If the kernel $f:\mathbb{R}^d \to \R_+$ is a positive definite
function that is positive on a set of non--zero Lebesgue measure
then $\kappa$ is positive definite.
\end{theorem}
\begin{proof}
By formula (\ref{eq:cov_integral}), we have
$$
\kappa(h)=\int_{\mathbb{R}^d} f(h-x) f^{\langle \alpha
-1\rangle}(-x) \, dx, \quad h\in \mathbb{R}^d.
$$
Thus for any $m \in \mathbb{N}$, $z_1,\ldots,z_m \in \mathbb{R}$,
$(z_1,\ldots,z_n)^\top \neq (0,\ldots,0)^\top$ and $s_1,\ldots,s_m
\in \mathbb{R}^d$ it holds
\begin{eqnarray*}
&&\sum_{i,j=1}^m \kappa(s_i-s_j) z_i z_j =  \int_{\mathbb{R}^d}
\sum_{i,j=1}^m f(s_i-s_j-x) z_i z_j f^{\langle \alpha
-1\rangle}(-x)\, dx>0.
\end{eqnarray*}
\end{proof}

An example of a process $X$ satisfying conditions of Theorem
\ref{th:COL_MA_pd} is the $S\alpha S$ Ornstein--Uhlenbeck process:
for any fixed $\lambda > 0$
$$X(t) = \int_\mathbb{R} e^{-\lambda(t-x)}\ind(t-x\geq 0)\, M(dx), \quad t \in \mathbb{R}.$$
By \cite[p.~138]{ST94},  we have
$\widehat{X}(t)=e^{-\lambda(t-t_n)} X(t_n)$ if $t_1< \ldots
<t_n<t$.

\begin{theorem}\label{th:COL_cont_MA}
If the covariation function $\kappa$ is positive definite and
continuous then the COL predictor is continuous.
\end{theorem}
\begin{proof}
Since $\kappa$ is positive definite, matrix $K$ is invertible, and
we have
\begin{equation*}
 \begin{pmatrix}
        \lambda_1 (t)\\
        \vdots \\
        \lambda_n (t)
          \end{pmatrix} = \begin{pmatrix}
 \kappa(0) & \cdots & \kappa(t_n-t_1) \\
 \vdots & \ddots & \vdots \\
 \kappa(t_n-t_1) & \cdots & \kappa(0)\end{pmatrix}^{-1}\begin{pmatrix}
                \kappa(t-t_1) \\
                \vdots \\
                \kappa(t-t_n)
\end{pmatrix}.
\end{equation*}
Since $\kappa$ is continuous, the weights $\lambda_1,\ldots,\lambda_n$ are continuous in $t$.
\end{proof}

\begin{exercise}
Show that continuous kernel functions with compact support yield a
continuous covariation function $\kappa$. Use the dominated
convergence theorem.
\end{exercise}

\subsubsection{The COL Predictor for Gaussian and sub--Gaussian Random Fields}

Let $X$ be a sub--Gaussian random field,  i.e.,
$X(t)=A^{1/2}G(t)$, $t\in\R^d$ where $A \sim
S_{\alpha/2}((\cos(\pi\alpha/4))^{2/\alpha},1,0)$ and
$\boldsymbol{G}$ is a zero mean stationary Gaussian field
independent of $A$.  In \cite[Example~2.7.4]{ST94}, it is shown
that for sub--Gaussian random fields, the covariation function is
given by
\begin{equation}
\kappa(h) = 2^{-\alpha/2} C(h)C(0)^{(\alpha-2)/2},\quad h\in\R^d,
\label{eq:covariation_function_sub_Gaussian}
\end{equation}
where $C(\cdot)$ is the covariance function of $G$.

It is easy to see that in this case the system (\ref{eq:COL_sys})
coincides with the simple kriging system \eqref{eq:system_kriging}
for $G$:
\begin{equation}\label{eq:COLSubGauss}
 \begin{pmatrix}
 C(0) & \cdots & C(t_n-t_1) \\
 \vdots & \ddots & \vdots \\
 C(t_n-t_1) & \cdots & C(0)\end{pmatrix} \begin{pmatrix}
        \lambda_1 \\
        \vdots \\
        \lambda_n
          \end{pmatrix} = \begin{pmatrix}
                C(t-t_1) \\
                \vdots \\
                C(t-t_n)
\end{pmatrix}.
\end{equation}
If $C$ is positive definite then the corresponding covariance
matrix is invertible which ensures the existence and uniqueness of
the solution of the system \eqref{eq:COLSubGauss}.

\begin{theorem}\label{th:COL_continuity}
If $(X(t_1), \ldots, X(t_n))^\top$ is full--dimensional and the
covariance function $C$ of the Gaussian component is continuous
then the COL predictor for sub--Gaussian random fields is
continuous.
\end{theorem}
The proof is similar to the proof of Theorem~\ref{th:COL_cont_MA}.

\begin{theorem}\label{th:LSL_COL}
Let $1<\alpha\leq 2$. For Gaussian and sub--Gaussian random
fields, the COL and LSL predictors coincide.
\end{theorem}
\begin{proof} Introduce the notation $t_0=t$. Put  $\lambda_0(t_0)=-1$ and
$$\widehat{X}(t_0)-X(t_0)=A^{1/2}\sum_{i=0}^n
\lambda_i(t_0)G(t_i).$$ The characteristic function of random
vector $(X(t_0),\ldots,X(t_n))^\top$ is given by
\begin{equation}\label{eq:charfctsubgaussian}
    \E \exp\left\{ i\sum_{k=0}^n \theta_k X(t_k)\right\} = \exp \left \{ -\left|\frac{1}{2}\sum_{i=0}^n \sum_{j=0}^n \theta_i \theta_j C (t_i-t_j)\right|
    ^{\alpha/2}\right\}
\end{equation}
for all $\theta_1,\ldots,\theta_n\in\R$, cf. \cite[Proposition
2.5.2]{ST94}. Now it is simple to see that
$$\sigma_{\widehat{X}(t_0)-X(t_0)}=\left(\frac{1}{2} \var\left(\sum_{i=0}^n \lambda_i(t_0)G(t_i)\right)\right)^{1/2}=
\left(\frac{1}{2} \sum_{i,j=0}^n\lambda_i\lambda_j C(t_i-t_j)
\right)^{1/2}.$$ Thus, the LSL optimization problem is equivalent
to
$$\sum_{i,j=0}^n\lambda_i\lambda_j C(t_i-t_j)\to \min_{\lambda_1,\ldots, \lambda_n}.$$
Taking derivatives we obtain $\sum_{j=0}^n C(t_k-t_j)
\lambda_j=0,$ $k= 1,\ldots, n$ which coincides with the COL
extrapolation system \eqref{eq:COLSubGauss}.
\end{proof}

\begin{remark}\label{rem:LSLsubGauss} It follows from the proof of Theorem \ref{th:LSL_COL} (which is valid for all $\alpha\in (0,2)$) that
the weights of the LSL predictor for sub--Gaussian random fields
are a solution of the system \eqref{eq:COLSubGauss} also in the
case $\alpha\in (0,1]$. The statement of Theorem
\ref{th:COL_continuity} holds as well. To summarize, the LSL
predictor for stationary sub--Gaussian random fields $X$ exists
and is unique and exact for all $\alpha\in (0,2]$ if the
covariance function $C$ of the Gaussian component $G$ is positive
definite. If $C$ is additionally continuous then this LSL
predictor is also continuous.
\end{remark}


\subsection{Maximization of Covariation}
\label{subsec:4.3} In this section, we assume that $X$ is an
$\alpha$--stable random field \eqref{eq:RFX} with
$\alpha\in(1,2]$. The predictor $\widehat{X}(t)=\sum_{i=1}^n
\lambda_i (t) X(t_i)$, whose weights
$\lambda_1(t),\ldots,\lambda_n(t)$  solve the
 following optimization problem
\begin{equation}
 \begin{cases}
  \left[\widehat{X}(t),X(t)\right]_\alpha= \sum_{i=1}^n \lambda_i (t)[X(t_i),X(t)]_\alpha \, \to \, \underset{\lambda_1,\ldots,\lambda_n}{\max}, \\
  \sigma_{\widehat{X}(t)} = \sigma_{X(t)}
 \end{cases}\label{MCL}
\end{equation}
 for $t\in \R^d$,  is called \textit{Maximization of
Covariation Linear (MCL) predictor}\index{Maximization of
Covariation Linear (MCL) predictor}.

The Lagrange function of the optimization problem \eqref{MCL} is
given by
\begin{equation*}
L(\vec{\lambda},\gamma) = \sum_{i=1}^n \lambda_i [X(t_i),X(t)]_\alpha + \gamma \left(\sigma^{\alpha}_{\sum_{i=1}^n \lambda_i X(t_i)} - \sigma^{\alpha}_{X(t)}\right), \ \vec{\lambda}\in\R^n,\ \gamma\in \R.
\end{equation*}
By taking partial derivatives and setting them equal to zero, we get
\begin{equation}
 \begin{cases}
  [X(t_j),X(t)]_\alpha + \gamma \cdot \partial \sigma^{\alpha}_{\sum_{i=1}^n \lambda_i X(t_i)}/\partial \lambda_j= 0, \quad j=1,\ldots,n, \\
  \sigma_{\sum_{i=1}^n \lambda_i X(t_i)} = \sigma_{X(t)}.
 \end{cases}\label{MCL_sys}
\end{equation}
Analogously to formula (\ref{eq:LSL_covariation}) one can obtain
$$\frac{\partial \sigma^{\alpha}_{\sum_{i=1}^n \lambda_i X(t_i)}}{\partial \lambda_j}= \alpha\cdot\left[X(t_j),\sum_{i=1}^n \lambda_i (t) X(t_i)\right]_\alpha.$$
Since $\gamma=-1/\alpha$, $\lambda_i(t_j)=\delta_{ij}$ is
obviously a solution of system (\ref{MCL_sys}) for $t=t_j$,
$j=1,\ldots, n$, the MCL predictor is exact.

Let us discuss the properties of the MCL predictor. Notice that
here no direct analogy with kriging can be drawn. For instance, a
counterpart $ \sigma_{\widehat{X}(t)} \le \sigma_{X(t)}$ of the
shrinkage property \eqref{eq:shrinkKr} is deliberately mutated to
the additional condition $ \sigma_{\widehat{X}(t)} =
\sigma_{X(t)}$. The reason for this is that both conditions lead
to the same solutions due to the convexity of the optimization
problem \eqref{MCL}.

Introduce the following  notation:
$\zeta(t)=(\kappa(t_1,t),\ldots, \kappa(t_n,t))^\top,$ $t\in\R^d,$
the function $\sigma_0: \R^d\to \R_+$ is $\sigma_0(t)=
\sigma_{X(t)}=\kappa(t,t)$. The function $\Psi: \mathbb{R}^n \to
\mathbb{R}_+$ is defined by
$$
\Psi(\lambda) = \sigma_{\widehat{X}(t)} = \left\Vert\sum_{i=1}^n \lambda_i f_{t_i} \right\Vert_\alpha.
$$
Denote the \emph{level set}\index{level set} of function $\Psi$ at
level $u\in\R$ by $ B_{u} = \{\lambda \in \mathbb{R}^n:
\Psi(\lambda) \leq u\}.$ The \emph{support set}\index{support set}
of any convex set $B\subset \R^n$ at a point $x\in \R^n$ is
defined by
$$
T(B,x)=\left\{ y \in B: \langle y,x \rangle = \sup _{z\in
B}\langle z, x \rangle\right\}.
$$
It is known that for strictly convex sets $B$ and any non--zero
$x\in \R^n$ the support set $T(B,x)$ is a singleton. We denote
this single point by $y_{B,x}$.

\begin{theorem}\label{th:MCL1}
Assume that the $\alpha$-stable random vector
$X=(X(t_1),\ldots,X(t_n))^\top$ is full--dimensional.
\begin{enumerate} \item The solution of the optimization problem
(\ref{MCL}) exists for all $t \in \mathbb{R}^d$. If
$\kappa(t_i,t)\neq 0$ for some $i =1,\ldots, n$ then the MCL
predictor $\widehat{X}(t)$ is unique. \item  If $\kappa$ is a
continuous function on $\mathbb{R}^d\times\mathbb{R}^d$ and
$\kappa(t_i,t)\neq 0$ for some $i =1,\ldots, n$ then the MCL
predictor is continuous in $t$.
\end{enumerate}
\end{theorem}
\begin{proof}\footnote{The authors are grateful to  D.Stolyarov and P.Zatitsky  for the help with the proof simplification.}
For the proof of the existence and uniqueness of MCL we refer the
reader to the paper \cite{KSS12}. It is also shown there that the
vector of MCL weights
$$\lambda(t)=(\lambda_1(t),\ldots,\lambda_n(t))^\top$$ is equal to
$y_{B_{\sigma_0(t)},\zeta(t)}$ for any $t\in \R^{d}$ whereas the
set $B_{\sigma_0(t)}$ is strictly convex. Let us prove that
$\lambda: \R^d\to \R^n$ is a continuous function. It is easy to
see that $B_{\sigma_0(t)}=\frac{1}{\sigma_0(t)}B_1$, because the
sets $B_{\sigma_0(t)}$, $t\in \R^d$ are homothetic, i.e. $a
B_{\sigma_0(t)}=B_{\sigma_0(t)/a}$, $a>0$. Thus by simple
geometric considerations
$$
T(B_{\sigma_0(t)}, \zeta(t))= T\left(\frac{1}{\sigma_0(t)}
B_1,\zeta(t)\right)=\frac{1}{\sigma_0(t)} T(B_1, \zeta(t)),
$$
thus $ \lambda(t)= \frac{1}{\sigma_0(t)}y_{B_1, \zeta(t)}$. Put
$B=B_1$ and $x(s)=y_{B, \zeta(s)}$ for any $s\in \R^d$. Show that
$\lim_{s\to t} x(s)= x(t)$. This limit exists by the definition of
the support set and continuity of the scalar product.  We know
that $\zeta(s)\to \zeta(t)$ as $s\to t$ since $\kappa$ is a
continuous function. Moreover, $B$ is a compact, and $x(s)\in B$
for all $s$. Choose a convergent sequence $s_m\to t$ as $m\to
\infty$ such that $x(s_m)\to y$ as $m\to \infty$, where $y\in B$.
Show that $y=x(t)$. It is clear that $\langle x(s_m),\zeta(s_m)
\rangle\to\langle y,\zeta(t)\rangle$ as $m\to \infty$. And for any
$x\in B$ it holds
$$
\langle x, \zeta(t)\rangle=\lim_{m\to\infty}\langle x,
\zeta(s_m)\rangle\leq \lim_{m\to\infty}\langle x(s_m),
\zeta(s_m)\rangle = \langle y, \zeta(t)\rangle.
$$
The inequality here is due to the fact that $\{x(s_m)\}=T(B,
\zeta(s_m))$ for any $m\in \N$. Thus $y=y_{B, \zeta(t)}$.
\end{proof}

\subsection{Case $\alpha\in(0,1]$}
\label{subsec:4.4}

As noticed in Section \ref{subsec:2.4}, the covariation function
is not defined for $\alpha \in (0,1]$. Moreover, the function
$\Vert \cdot \Vert_\alpha$ for $\alpha<1$ defined in
\eqref{eq:L_alpha_norm} is not a norm anymore since the triangle
inequality fails to hold. The property of strict convexity of
$L^\alpha(E,m)$ does not hold as well.

To cope with these drawbacks, one may come to an idea that the
codifference (cf. Definition \ref{dfn:codiff}) can be used instead
of the covariation in COL and MCL methods.  However, it does not
seem to make advances in extrapolation. For instance, replacing
the covariation by the codifference in the MCL method leads to the
 optimization problem
\begin{equation}
 \begin{cases}\label{eq:MCL01}
  \tau (\widehat{X}(t), X(t)) = \sigma_{\widehat{X}(t)} + \sigma_{X(t)} - \sigma_{\widehat{X}(t)-X(t)}
   \ \to \ \max\limits_{\lambda_1,\ldots,\lambda_n} \\
  \sigma_{\widehat{X}(t)} = \sigma_{X(t)}.
 \end{cases}
\end{equation}
Using the constraint $\sigma_{\widehat{X}(t)} = \sigma_{X(t)}$,
the first relation rewrites
\begin{equation*}
    \tau (\widehat{X}(t), X(t)) = 2\sigma_{X(t)} - \sigma_{\widehat{X}(t)-X(t)}.
\end{equation*}
Hence, the method \eqref{eq:MCL01} is equivalent to LSL
extrapolation, i.e., to minimizing the scale parameter
\begin{equation*}
    \sigma_{\widehat{X}(t)-X(t)} = \|f_{t}-\sum_{i=1}^n\lambda_i
     f_{t_i}\|_{\alpha}
\end{equation*}
of $\widehat{X}(t)-X(t)$.

Replacing  the covariation by the codifference in the COL
    method \eqref{eq:COL},
    one arrives at the system of nonlinear equations
    \begin{equation}\label{eq:orthogonalCodiff}
        \tau_{\widehat{X}(t), X(t_i)} = \tau_{X(t), X(t_i)} , \ \ \ i=1,\ldots,n.
    \end{equation}
    Here the numerical computation of a solution
    is necessary, which can be very time consuming. Furthermore, it is shown in \cite{Hagel12} that the solution of the system \eqref{eq:orthogonalCodiff} is
    not unique.
    For this reason, we shall not pursue the method (\ref{eq:orthogonalCodiff}) in
    future.

    Neither leads the maximization of $\tau_{\widehat{X}(t),
    X(t)}$ with respect to weights $\lambda_1,\ldots,\lambda_n$ to a unique
    predictor \eqref{eq:LinPred}. In particular, its existence is not really clear.
    As an example consider a random field \eqref{eq:RFX} with the
    kernel function $f_t$ of compact support such that the
    supports of $f_t$ and $f_{t_1},\ldots, f_{t_n}$ do not
    overlap. Then it is easy to see that $\tau_{\widehat{X}(t),
    X(t)}=0$ allowing for an arbitrary choice of weights
    $\lambda_1,\ldots,\lambda_n$.

In the remainder of this Section, we focus on the properties of
the LSL method for $\alpha$--stable random fields with $\alpha \in
(0,1]$. First of all, the fundamental question of existence has to
be answered. Here we follow \cite{Hagel12} and do this in a more
general setting of $r$--normed vector spaces.
\begin{definition}
Let $V$ be a vector space over a field $\mathbb{K}$. A map
$||.||_{(r)}:V \rightarrow \mathbb{R}_+$  is called an
\emph{$r$--norm}\index{$r$--norm}, if there exists $K \geq 1$ and
$r > 0$ such that
\begin{align*}
     ||x||_{(r)} & =  0 \Leftrightarrow x = 0, \\
     ||ax||_{(r)} & =  |a| \cdot ||x||_{(r)} \ \ \ \forall \ a \in \ \mathbb{K}, \ \forall \ x \ \in \ V, \\
     ||x+y||_{(r)} & \leq K(||x||_{(r)}+||y||_{(r)}) \ \ \ \forall \ x,y \ \in \ V, \\
     ||x+y||_{(r)}^r & \leq ||x||_{(r)}^r + ||y||_{(r)}^r \ \ \ \forall \ x,y \ \in \ V.
\end{align*}
\end{definition}
Now the existence theorem can be formulated.
\begin{theorem}[\cite{Hagel12}] \label{theo:LSLexist}
Let $V$ be a vector space over $\mathbb{R}$ with r-norm $|| \cdot
||_{(r)}$ and let $f_1,\ldots,f_n \in V$ be linearly independent.
For any $f_0\in V$, there exist real numbers
$\lambda_1^*,\ldots,\lambda_n^*$ such that
\begin{equation*}
    ||f_0-\sum_{i=1}^n \lambda_i^{*}f_i||_{(r)} = \inf\limits_{\lambda_1,\ldots,\lambda_n \in \mathbb{R}}
    ||f_0-\sum_{i=1}^n \lambda_i f_i||_{(r)}.
\end{equation*}
\end{theorem}
If we set $V = L^{\alpha}(E,m)$ and note that $|| \cdot
||_{(\alpha)} = \Vert \cdot \Vert_\alpha$ defined in
\eqref{eq:L_alpha_norm} is an $\alpha$-norm on $L^{\alpha}(E,m)$
(even a norm if $\alpha \geq 1$), the existence of the LSL
predictor follows immediately from Theorem \ref{theo:LSLexist}. In
contrast to the case $\alpha \in (1,2]$ (Theorem \ref{th:LSL1}),
the uniqueness of the LSL weights $\lambda^* :=
(\lambda_1^*,\ldots,\lambda_n^*)^\top$ in Theorem
\ref{theo:LSLexist} is not guaranteed.  We illustrate this by the
following example. Introduce the notation
$H_{\alpha}(\lambda)=\sigma_{\widehat{X}(t)-X(t)}$ for
$\lambda=(\lambda_1,\ldots,\lambda_n)^\top \in \R^n$.
\begin{example}
Consider the measurable space $(E, m) = ([0,1], \nu_1)$ and the
kernel function $f_t(x) = \ind
\left(x\in(t+\frac{1}{4},t+\frac{3}{4})\right)$. Given $t_1 =
\frac{1}{4}$, predict the value of the symmetric $\alpha$--stable
process $X(t) = \int\limits_{[0,1]}f_t(x) \, M(dx)$ at the point
$t = 0$. By elementary calculations we obtain
\begin{equation*}
    H_{\alpha}^{\alpha}(\lambda) = \int\limits_{[0,1]}|f_{t}(x)-\lambda f_{t_1}(x)|^{\alpha} dx = \frac{1}{4}(1+|1-\lambda|
    ^{\alpha} + |\lambda|^{\alpha}).
\end{equation*}
It is easy to see that for $0<\alpha < 1$, $H_{\alpha}$ has two
global minima at $\lambda=0$ and $\lambda=1$. If $\alpha=1$ the
set of all global minimum points equals the interval $[0,1]$. For
values $\alpha
> 1$, the function $H_{\alpha}$ has a unique global
minimum at $\lambda = 0.5$.
\end{example}

In order to get unbiased prediction (provided that the first
moment of $X$ is finite), the parameter space is often restricted
to $\{(\lambda_1,\ldots ,\lambda_n)^\top \in \R^n:  \ \
\sum_{i=1}^{n}\lambda_i = 1 \}$. S. Hagel showed in \cite{Hagel12}
that this restriction does not cause uniqueness of LSL prediction
for $\alpha\in (0,1)$. Alternatively, the following algorithmic
approach to choose a unique global minimum in the LSL optimization
problem is proposed:
\begin{algo}\label{algo:bestLSL}
    Let $\{X(t),\; t \in T\}$ be an $\alpha$--stable  random field \eqref{eq:RFX}
    with $0 < \alpha < 1$ and
    $T \subset \R ^d$. Let $t_1,\ldots,t_n \in T$ be fixed such that functions
    $f_{t_1},\ldots,f_{t_n}$ are linearly independent.
    \begin{enumerate}
    \item Order the points $t_1,\ldots,t_n$ so that
    \begin{equation*}
        \|t-t_1\| \leq \|t-t_2\| \leq \ldots \leq \|t-t_n\|
    \end{equation*}
    and if $\|t-t_i\| = \|t-t_{i+1}\|$ for some $i \in \{1,\ldots,n-1\}$ then
    \begin{align}
         t_i^{(p)}  & =t_{i+1}^{(p)} \ \ \ \ \mbox{for all } p=1,\ldots,k-1 \label{eq:unicond1} \\
         t_i^{(k)}  & < t_{i+1}^{(k)} \label{eq:unicond2}
    \end{align}
    for some $k \in \{1,\ldots,m\}$, where $t_i^{(p)}$ is the $p$--th component of $t_i$.
    \item Determine the set $A_0$ of all critical points
    \begin{equation*}
        A_0 = \{(\lambda_1,\ldots,\lambda_n) \in \R^n \ \text{:} \ \ H_{\alpha}(\lambda_1,\ldots,\lambda_n) = \inf\limits_
        {(\mu_1,\ldots,\mu_n) \in \R^n} H_{\alpha}(\mu_1,\ldots,\mu_n)\}
    \end{equation*}
    \item Reduce $A_0$ step by step to sets $A_1 \supseteq A_2 \supseteq \dots \supseteq A_n$ given by
    \begin{equation*}
        A_j = \{(\lambda_1,\ldots,\lambda_n) \in A_{j-1} \ \text{:} \ \ \lambda_j = \max\limits_
        {(\mu_1,\ldots,\mu_n) \in A_{j-1}} \mu_j\}, \ \ \ \  \ j=1,\ldots,n.
    \end{equation*}
    \end{enumerate}
    Clearly, the set $A_n$ consists of just one element.
\end{algo}

\begin{definition}\label{def:bestLSL}
    We call
    $ \widehat{X}(t) = \sum_{i=1}^n \lambda_i^* X(t_i)$
    the \emph{best LSL predictor}\index{best LSL predictor} if $(\lambda_1^*,\ldots,\lambda_n^*) \in A_n$.
\end{definition}

The above construction has a simple intuitive meaning. The points
$t_1,\ldots,t_n$ are ordered with respect to their distance to
$t$. To get a unique ordering, conditions (\ref{eq:unicond1}) and
(\ref{eq:unicond2}) are required. Points with a smaller distance
to $t$ are regarded to exert more influence on the value of $X$ at $t$, so their weights should be maximized first. \\

To show that $A_j\neq \emptyset$, $j = 1,\dots,n$ we notice that
$A_0$ is nonempty and compact. Therefore, the projection mapping
$(x_1,\ldots,x_n) \mapsto x_1$ takes its maximum on $A_0$. Hence,
$A_1$ is nonempty and compact as well. Sets $A_2,\ldots,A_n$ are
not empty by induction.

It can be easily proved that the best LSL predictor is exact. To
see this, let $t = t_i$ for some $i \in \{1,\ldots,n\}$ and let
$t_1,\ldots,t_n \in \R^d$ be as in Algorithm \ref{algo:bestLSL}.
    Relations (\ref{eq:unicond1}) and (\ref{eq:unicond2}) then imply that $t = t_1$. Trivially
    $(1,0,\ldots,0) \in A_0$ holds. Due to the linear independence of $f_{t_1},\ldots,f_{t_n}$, it holds that
    $A_n=A_0 = \{(1,0,\ldots,0)\}$.

For $1 < \alpha \leq 2$,  Theorem \ref{th:LSL2} stated the
continuity of LSL prediction. In contrast, the best LSL predictor
is not necessarily continuous for $0 < \alpha \le 1$ as the next
example shows.
\begin{example}
    Let $X=\{X(t), \; t\in\R^2\} $ be an $\alpha$--stable random field \eqref{eq:RFX} with  $0 < \alpha < 1$,
$f_{t}(x) =\ind \left(x\in (\min\{t^{(1)},t^{(2)}\},
\max\{t^{(1)},t^{(2)}\})\right)$ for $t = (t^{(1)},t^{(2)}) \in
\R^2$, $E=\R$ and $M$ being a $S\alpha S$ random measure on $\R$
with Lebesgue control measure. It follows from relations
\eqref{eq:p_mean}, \eqref{eq:scale} and Markov inequality that
 that $X$ is stochastically continuous, i.e., it has
a.s. no jumps at fixed locations $t$. For
    $n = 1$, introduce $t_0 = (\frac{1}{2}, \frac{3}{2})$, $t_1 = (0, 1)$, $t=t_0+\varepsilon$, where  $\varepsilon
    = (\delta, \delta) \in \R^2$ for some $\delta \in (-\frac{1}{2}, \frac{1}{2})$.
    Consider the best LSL predictor $\widehat{X}(t)$ of $X(t)$ based on the data
    $X(t_1)$. It holds
    \begin{align*}
        H_{\alpha}^{\alpha}(\lambda) & = \int\limits_{\R} |f_{t_0+\epsilon}(x)-\lambda f_{t_1}(x)|^{\alpha} dx \\
                            & = \left(\frac{1}{2} + \delta \right) \cdot |\lambda|^{\alpha} +  \left(\frac{1}{2} - \delta \right)
                                 \cdot |1-\lambda|^{\alpha} +  \left(\frac{1}{2} + \delta \right).
    \end{align*}
    If $\delta > 0$ then $H_{\alpha}$ has a global minimum at $\lambda = 0$ and if $\delta < 0$ it has a
    global minimum at $\lambda = 1$. So $\widehat{X}(t)$ is discontinuous at $t=t_0$.
\end{example}


In addition to the best LSL prediction, it is possible to treat
the case $\alpha = 1$ similar to the case $1 < \alpha < 2$. The
following approach is proposed in \cite{Hagel12}. For a symmetric
$1$--stable field $\{X(t) \ \text{:} \ t \in T\}$ with integral
representation
\begin{equation*}
    X(t) = \int\limits_{E} f_t(x)\, M(dx)
\end{equation*}
let the function $f_t \in L^1(E,\mathcal{E},m) \cap
L^{\delta}(E,\mathcal{E},m)$ for some $\delta > 1$. Then we have
\begin{equation*}
    \int\limits_{E} |f_{t}(x)-\sum_{i=1}^n \lambda_i f_{t_i}(x)|^{\gamma}m(dx) < \infty
\end{equation*}
for all $\gamma \in [1,\delta]$, $\lambda_1,\ldots,\lambda_n \in
\R$ and $t,t_1,\ldots,t_n \in T$. Now fix $t,t_1,\ldots,t_n \in T$
and chose an arbitrary sequence $(\gamma_k)_{k \in \N} \subset
(1,\delta]$  which converges to $1$ as $k\to\infty$. Let
$(\lambda_1^{(\gamma_k)},\ldots,\lambda_n^{(\gamma_k)})$ be the
unique solution of
\begin{equation}\label{eq:LSL_gamma_k}
    \int\limits_{E} |f_{t}(x)-\sum_{i=1}^n \lambda_i f_{t_i}(x)|^{\gamma_k}m(dx) \rightarrow \min_{\lambda_1,\ldots,
    \lambda_n}.
\end{equation}
Applying the stability theorem in \cite[p.225]{Kos91} it follows
the convergence
\begin{equation}\label{eq:convStability}
    \int\limits_{E} |f_{t}(x)-\sum_{i=1}^n \lambda_i^{(\gamma_k)} f_{t_i}(x)|^{\gamma_k}m(dx) \rightarrow \inf_{\mu_1,\ldots,
    \mu_n} \int\limits_{E} |f_{t}(x)-\sum_{i=1}^n \mu_i f_{t_i}(x)|m(dx)
\end{equation}
as $k\to\infty.$ Moreover, it can be shown that
\begin{equation*}
    (\lambda_1^{(\gamma_k)},\ldots,\lambda_n^{(\gamma_k)}) \rightarrow
    (\lambda_1^*,\ldots,\lambda_n^*),
    \quad  k \rightarrow \infty.
\end{equation*}
This set of weights $(\lambda_1^*,\ldots,\lambda_n^*)$ exists and
is unique\footnote{Personal communication of Adrian Zimmer} if all
LSL prediction problems \eqref{eq:LSL_gamma_k} with stability
indices $\gamma_k>1$  do so. It also does not depend on the choice
of the sequence $(\gamma_k)_{k \in \N} \subset (1,\delta]$ such
that $\gamma_k \to 1$ as $k\to\infty$.
\begin{definition}
The predictor $\widehat{X}^*(t)=\sum_{i=1}^n \lambda_i^* X(t_i)$,
$t\in T$ is called  an \emph{index--continuous LSL predictor}
(ICLSL)\index{index--continuous LSL predictor} for the symmetric
$1$--stable random field $X$.
\end{definition}
It is still an open problem to explore the statistical properties
of ICLSL.

\subsection{Numerical Examples}
\label{subsec:4.5}

In this section, LSL, COL and MCL extrapolation methods (as well
as Maximum Likelihood extrapolation and conditional simulation for
sub--Gaussian random fields) are applied to simulated data of
various $\alpha$--stable random processes and fields $X$ for
$\alpha\in (0,2)$.

The random fields are simulated and extrapolated on an equidistant
$50 \times 50$ --grid of points within $T=[0,1]^2$. In Examples 1
and 2, the simulated field $X=\{  X(t),\; t\in [0,1]^2\}$ is
observed at the points $t_1,\ldots,t_{16}$ given by their
coordinates
\begin{align*}
     t_1 & = (0,0), & t_2 & = (0,0.3), & t_3 & = (0,0.6), & t_4 & = (0,0.9), \\
     t_5 & = (0.3,0), & t_6 & = (0.3,0.3), & t_7 & = (0.3,0.6), & t_8 & = (0.3,0.9), \\
     t_9 & = (0.6,0), & t_{10} & = (0.6,0.3), & t_{11} & = (0.6,0.6), & t_{12} & = (0.6,0.9), \\
     t_{13} & = (0.9,0), & t_{14} & = (0.9,0.3), & t_{15} & = (0.9,0.6), & t_{16} & = (0.9,0.9).
\end{align*}

\begin{figure}[t!]

\subfigure[Realization of a sub--Gaussian random field with
$\alpha=1.2$.]{ \label{fig:100}
\includegraphics[scale=.30]{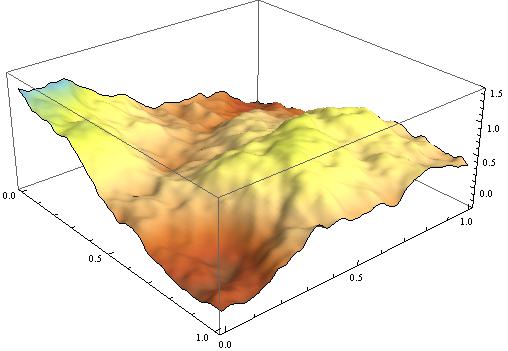}
} \hfill \subfigure[Corresponding LSL (COL) predictor]{
\label{fig:12}
\includegraphics[scale=.30]{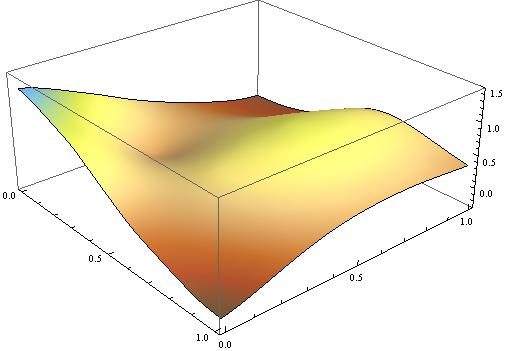}
} \subfigure[Corresponding MCL predictor]{ \label{fig:13}
\includegraphics[scale=.30]{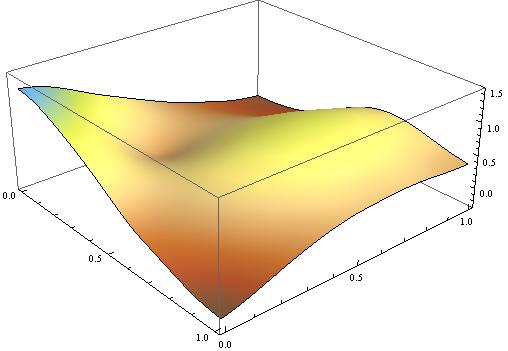}
} \hfill \subfigure[Pointwise difference ((b)-(c)) between LSL (b)
and MCL (c) predictors ]{ \label{fig:DiffLSLMCL}
\includegraphics[scale=.30]{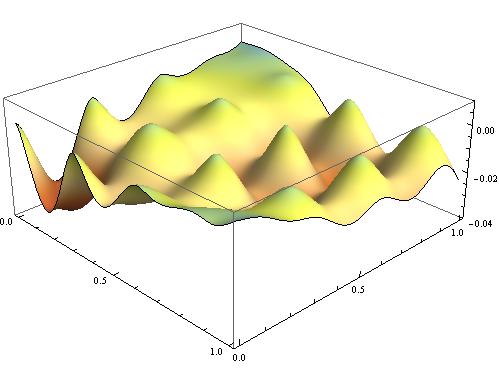}
} \caption{Realization of a sub--Gaussian random field for $\alpha
> 1$ and different predictors}

\end{figure}

\begin{figure}[t!]
\subfigure[Realization of a sub--Gaussian random field with
$\alpha = 0.8$.]{ \label{fig:14}
\includegraphics[scale=.30]{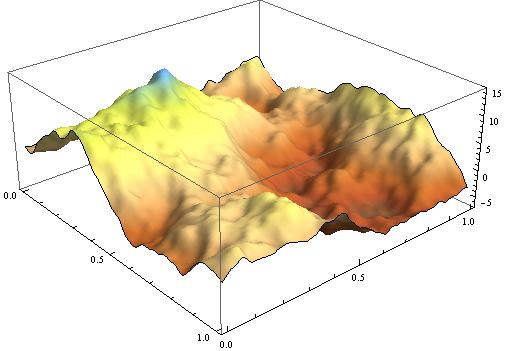}
} \hfill \subfigure[Corresponding LSL (ML) predictor]{
\label{fig:15}
\includegraphics[scale=.30]{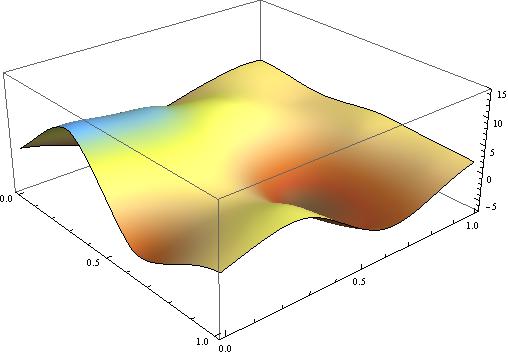}
} \centering \subfigure[Prediction by conditional simulation]{
\label{fig:16}
\includegraphics[scale=.30]{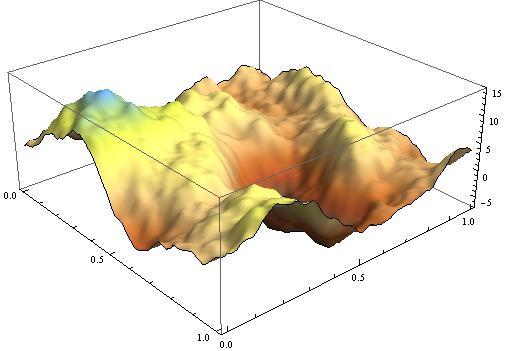}
} \caption{Realization of a sub--Gaussian random field for $\alpha
< 1$ and different predictors}
\end{figure}

\paragraph{\bf 1. Sub--Gaussian Random Fields}
Consider a stationary sub--Gaussian random field $X$ described in
Example 5 of Section \ref{subsec:2.5} with $\alpha = 1.2$. The
Gaussian part $G$ of this field has a Whittle--Mat\'ern covariance
function (cf. Section \ref{subsec:2.2.1}, Example 6) with
parameters as in Figure \ref{fig:2}. Figure \ref{fig:100} shows a
realization of $X$. The corresponding LSL (coinciding with COL by
Theorem \ref{th:LSL_COL}) and MCL predictors can be seen in
Figures \ref{fig:12} and \ref{fig:13}. Both predictions are
smoother than the realization of the field itself. Since
predictions in Figures \ref{fig:12} and \ref{fig:13} look quite
similar and can not be told one from another by eye, their
difference is given in Figure \ref{fig:DiffLSLMCL}.
 \\
Figure \ref{fig:14} shows a realization of the stationary
sub--Gaussian field with $\alpha = 0.8$ and covariance function
$C$ of the Gaussian part as above. A Maximum Likelihood (ML)
predictor for sub--Gaussian random fields is introduced in
\cite{KSS12}. It is shown in Theorem 11 of that paper that LSL,
COL and ML methods coincide if $\alpha\in (1,2)$. However, its
proof does not depend on $\alpha$ covering (with regard to Remark
\ref{rem:LSLsubGauss} of this chapter) the range of all
$\alpha\in(0,2)$. Thus, LSL and ML predictors coincide for
sub--Gaussian random fields with any stability index
$\alpha\in(0,2)$. A possibility of extrapolation of sub--Gaussian
random fields $X$ by conditional simulation (CS) of the Gaussian
component $G$ of $X$ and the subsequent scaling by $\sqrt{A}$ is
straightforward; see e.g.  \cite{Pai98} and
\cite[p.~112]{Karcher12}. Algorithms for the conditional
simulation of $G$ are given in \cite{Lan02}. Corresponding
extrapolation results for LSL (ML) and CS methods are given in
Figures \ref{fig:15} and \ref{fig:16}. Notice that the ML
prediction  for this realization of $X$ is much smoother than CS
prediction.

\begin{figure}[t!]
\subfigure[Realization of stable L\'{e}vy motion with skewness
intensity $\beta = 1$  and $\alpha = 1.5$]{ \label{fig:17}
\includegraphics[scale=.28]{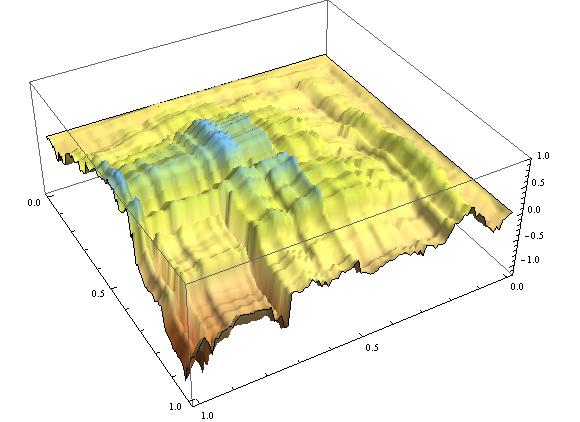}
} \hfill \subfigure[Corresponding LSL predictor]{ \label{fig:18}
\includegraphics[scale=.28]{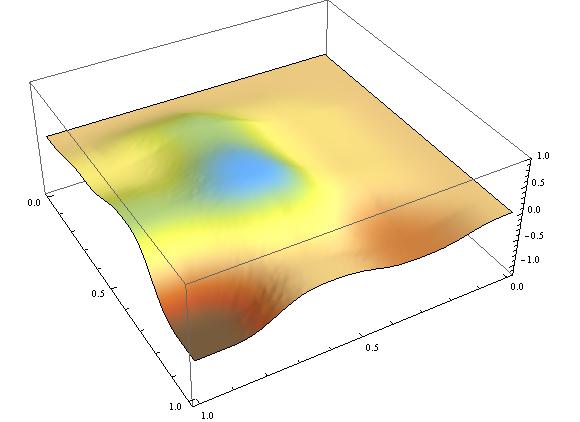}
} \subfigure[Corresponding COL predictor]{ \label{fig:19}
\includegraphics[scale=.28]{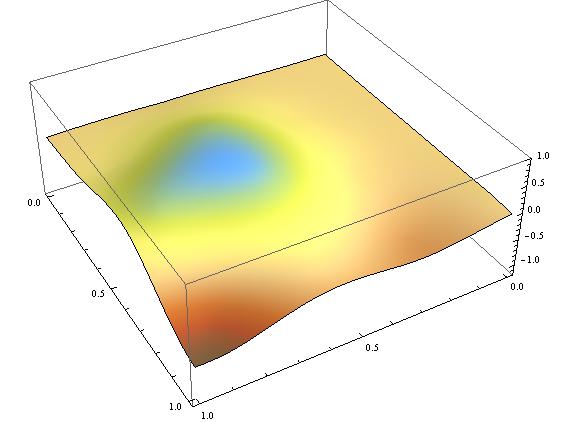}
} \hfill \subfigure[Corresponding MCL predictor]{ \label{fig:20}
\includegraphics[scale=.28]{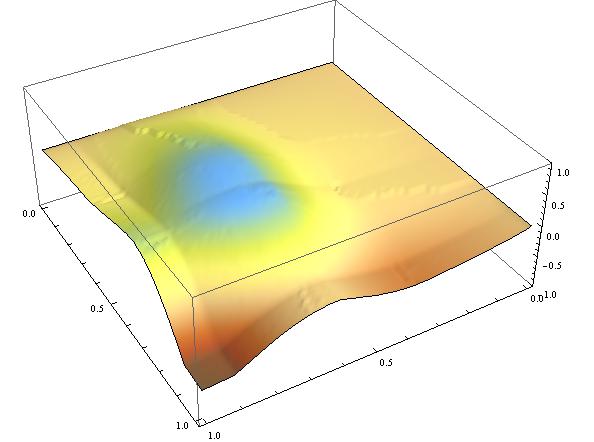}
} \caption{Realization of a skewed stable L\'{e}vy motion field
and different predictors}
\end{figure}

\paragraph{\bf 2. Skewed stable L\'{e}vy Motion}
Consider the two--dimensional $1.5$--stable  L\'evy motion $X$
defined by
$$X(t) = \int_{0}^1 \int_{0}^1 \ind (x_1\leq t_1,x_2 \leq
t_2)\,  M\big(d(x_1,x_2)\big),\quad t=(t_1,t_2)^\top \in
[0,1]^2,$$ where $M$ is a non--symmetric centered $1.5$--stable
random measure with skewness intensity $\beta = 1$. Comparing a
realization of $X$ (Figure \ref{fig:17}) with its LSL, COL and MCL
predictors (Figures \ref{fig:18}, \ref{fig:19} and \ref{fig:20})
one can see that prediction has a smoothing effect.

\begin{figure}[t!]
\includegraphics[scale=.5]{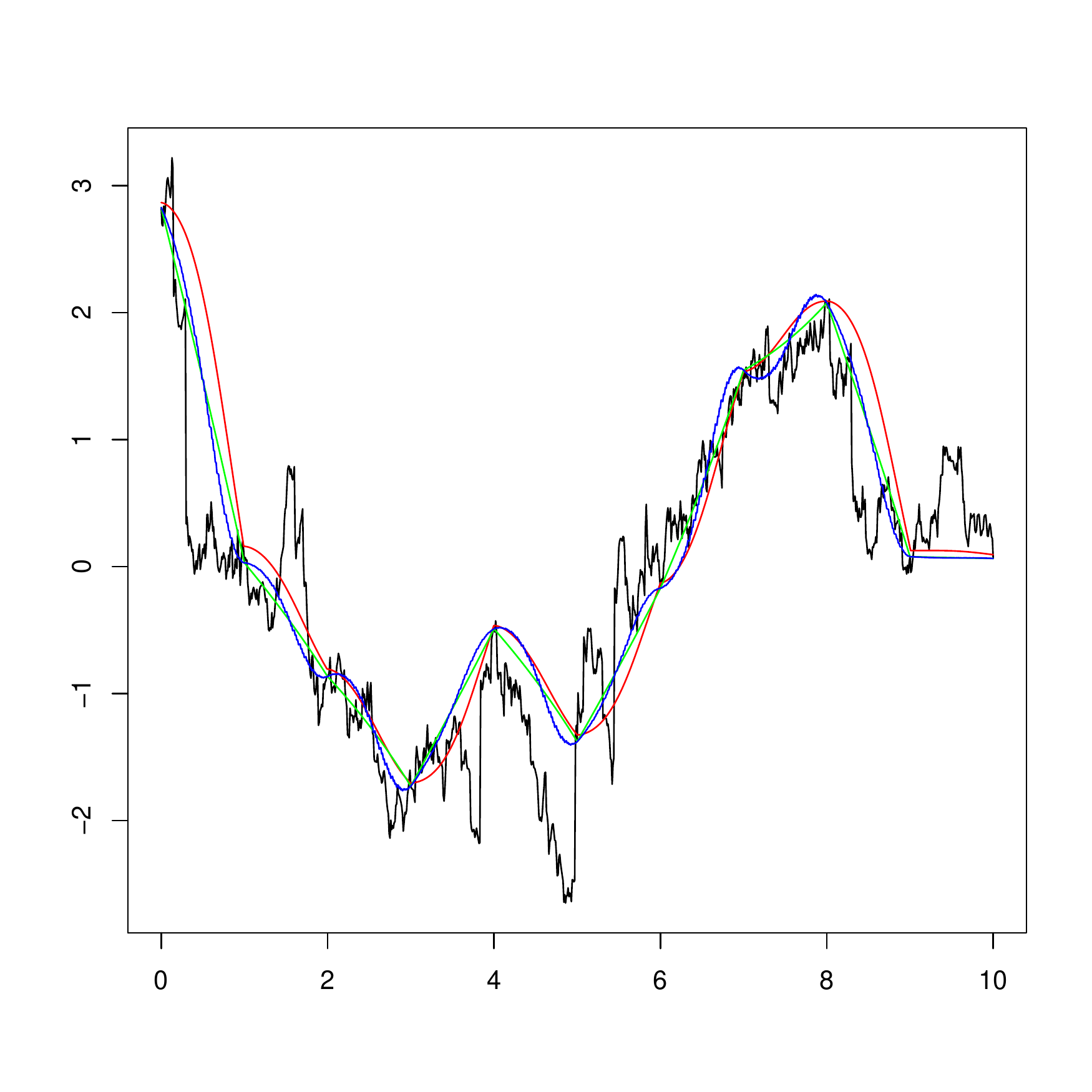}
\sidecaption \caption{A trajectory (black) of the stable
Ornstein--Uhlenbeck process together with LSL (red), COL (green)
and MCL (blue) predictors, $\alpha=1.6$} \label{fig:21}
\end{figure}

\paragraph{\bf 3. Stable Ornstein--Uhlenbeck Process}
Let $X$ be a $1.6$-stable Ornstein--Uhlenbeck process with
$\lambda=0.5$ defined in Example 2 of Section \ref{subsec:2.5}.
Figure \ref{fig:21} shows a trajectory of this process and
different interpolators. The process $X$ is observed at positions
$t_i=1,\dots,10$ within $[0,10]$. It can be seen that LSL
interpolation is very smooth. In contrast, the COL predictor is
piecewise smooth and continuous on the whole interval.

\begin{figure}[t!]
\subfigure[Realization of a skewed $0.5$--stable moving average
random field]{ \label{fig:SaSMA}
\includegraphics[scale=.3]{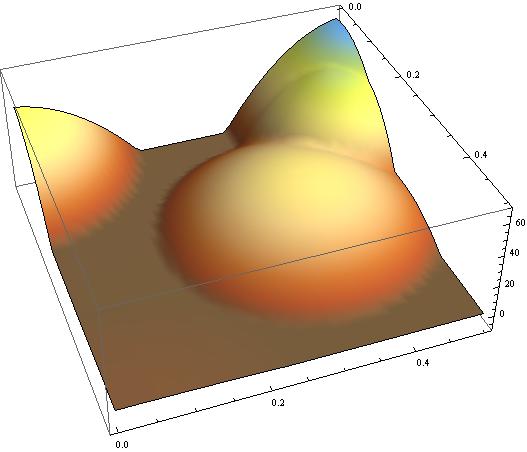}
} \hfill \subfigure[Corresponding best LSL predictor]{
\label{fig:SaSMAbestLSL}
\includegraphics[scale=.3]{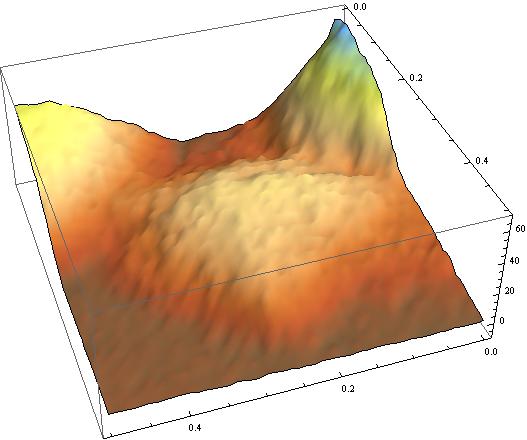}
}
 \caption{Realization of a skewed moving average
    field with $\alpha=0.5$ and its best LSL
    predictor}
\end{figure}

\paragraph{\bf 4. Stable Moving Average}
    Let $X=\{X(t), \ t \in [0,0.49]^2\}$ be a moving average
    field (cf. Example 2 of Section \ref{subsec:2.5}) with the
    kernel function
    \begin{equation*}
        f(x) = 0.5\left(0.04-\|x\|^2\right)\ind \left( \|x\| \leq
        0.2\right),
    \end{equation*}
     stability index $\alpha=0.5$ and skewness intensity $\beta = 0.8$.
    Random field $X$ is simulated on an equidistant $50 \times 50$--grid of points within $[0,0.49]^2$ using the step function approach from paper
    \cite{KSS13} with an accuracy ($L^{\alpha}$-error) $\epsilon = 0.01$. The field is observed
    at points
    \begin{align*}
     t_1 & = (0,0), & t_2 & = (0,0.25), & t_3 & = (0,0.49), \\
     t_4 & = (0.25,0), & t_5 & = (0.25,0.25), & t_6 & = (0.25,0.49), \\
     t_7 & = (0.49,0), & t_8 & = (0.49,0.25), & t_9 & =
     (0.49,0.49).
    \end{align*}
    To solve the optimization problems for the best LSL prediction (cf. Section
    \ref{subsec:4.4}) numerically, an average of $8$ realizations of the simulated annealing algorithm from \cite{Kirkpatrick1983} is used. Figures \ref{fig:SaSMA} and \ref{fig:SaSMAbestLSL} show  a realization of $X$
    and its best LSL predictor. The numerical optimization procedure is
    quite time consuming with $136$ min. of computation time  (Pentium Dual Core E5400, $2.70$ GHz, $8$ GB RAM) per
    extrapolation.



\section{Open problems}
\label{sec:5} In contrast to kriging methods, there is no common
methodology of measuring prediction errors in the stable case. We
propose the following measures
\begin{equation}\label{eq:openpr1}\sup_{t\in \mathbf{R}^d}\left( \mathbf{E}|X(t)-\widehat{X}(t)|^p \right)^{1/p}=
 c_\alpha (p) \sup_{t\in \mathbf{R}^d} \|f_{t}-\sum_{i=1}^n \lambda_i f_{t_i}\|_{\alpha},\end{equation}
where $1<p<\alpha$ and $c_\alpha (p)>0$ is a constant from
relation \eqref{eq:p_mean}, or
\begin{equation}\label{eq:openpr2}\mathbf{P}\left(\sup_{t\in \mathbf{R}^d} |X(t)-\widehat{X}(t)|
> \epsilon\right), \quad \epsilon>0.  \end{equation}
It is an open problem to find lower and upper bounds for these
errors as well as minimax bounds where the infimum over a subclass
of stable random fields $X$ is additionally considered in
relations \eqref{eq:openpr1} and \eqref{eq:openpr2}.
Alternatively, one can be interested in the asymptotic behavior of
$ \mathbf{P}\left(\sup_{t\in \mathbf{R}^d} |X(t)-\widehat{X}(t)| <
\epsilon\right)$ as $\epsilon\to 0$ which is related to small
deviation problems.

\begin{acknowledgement}
This research was partially supported by the DFG -- RFBR grant
09--01--91331. The second author was also supported by the
Chebyshev Laboratory (Department of Mathematics and Mechanics,
St.-Petersburg State University) within RF government grant
11.G34.31.0026.
\end{acknowledgement}


\bibliographystyle{spmpsci}
\bibliography{ssr_bibl}

\printindex

\end{document}